\documentclass[11pt]{amsart}

\usepackage{amsmath}
\usepackage{amssymb}

\newtheorem{theorem}{Theorem}[section]
\newtheorem{claim}[theorem]{Claim}

\newtheorem{lemma}[theorem]{Lemma}

\newtheorem{proposition}[theorem]{Proposition}
\newtheorem{corollary}[theorem]{Corollary}

\theoremstyle{definition}
\newtheorem{definition}[theorem]{Definition}

\newtheorem{question}[theorem]{Question}

\theoremstyle{remark}
\newtheorem{remark}[theorem]{Remark}

\newcount\skewfactor
\def\mathunderaccent#1#2 {\let\theaccent#1\skewfactor#2
\mathpalette\putaccentunder}
\def\putaccentunder#1#2{\oalign{$#1#2$\crcr\hidewidth
\vbox to.2ex{\hbox{$#1\skew\skewfactor\theaccent{}$}\vss}\hidewidth}}
\def\name{\mathunderaccent\tilde-3 }

\def\smallbox#1{\leavevmode\thinspace\hbox{\vrule\vtop{\vbox
   {\hrule\kern1pt\hbox{\vphantom{\tt/}\thinspace{\tt#1}\thinspace}}
   \kern1pt\hrule}\vrule}\thinspace}

\newcommand{\cf}{{\rm cf}}

\DeclareMathOperator{\crit}{crit}

\def\qedref#1{$\qed_{\reforiginal{#1}}$}

\title{Magidor Cardinals}
\author{Shimon Garti}
\address{Institute of Mathematics,
 The Hebrew University of Jerusalem,
 Jerusalem 91904, Israel}
\email{shimon.garty@mail.huji.ac.il}

\author{Yair Hayut}
\address{Institute of Mathematics,
 The Hebrew University of Jerusalem,
 Jerusalem 91904, Israel}
\email{yair.hayut@mail.huji.ac.il}


\begin{document}
\let\labeloriginal\label
\let\reforiginal\ref

\begin{abstract}
We define Magidor cardinals as J\'onsson cardinals upon replacing colorings of finite subsets by colorings of $\aleph_0$-bounded subsets. Unlike J\'onsson cardinals which appear at some low level of large cardinals, we prove the consistency of having quite large cardinals along with the fact that no Magidor cardinal exists.
\end{abstract}

\maketitle

\newpage

\section{Introduction}
The definitions of large cardinals came, historically, from different fields in mathematics. Some notions are set-theoretical (e.g., inaccessible or Mahlo cardinals), some belong to the logic realm (e.g., compact and weakly compact cardinals). Many of them were defined on pure combinatorial grounds (e.g., Ramsey or Erd\H{o}s cardinals) and one definition comes from measure theory (i.e., the measurable cardinals). For most of large cardinals, a natural defining property exists in more than one field.

The higher part of the chart of large cardinals is connected with elementary embeddings. If $\kappa$ is large enough then there is a non-principal ultrafilter $U$ on $\kappa$ which is $\kappa$-complete, and one can form the ultraproduct ${\rm V}^\kappa/U$ and get a well-founded model. The construction yields an elementary embedding $\jmath:{\rm V}\rightarrow {\rm V}^\kappa/U$ which is not the identity.

An important question here is how far can we stretch the similarity between the original universe ${\rm V}$ and the resulting model ${\rm V}^\kappa/U$.
More generally, if $\jmath:{\rm V}\rightarrow M$ and $\kappa=\crit(\jmath)$ then $\kappa$ is a large cardinal (regardless if it comes from an ultraproduct construction), and we can ask how far is ${\rm V}$ from $M$.
It turns out that this question produces stronger and stronger notions of large cardinals (the main directions being strong cardinals on one hand, and supercompact or huge cardinals on the other hand). The ultimate demand is ${\rm V}={\rm V}^\kappa/U$ (or ${\rm V}=M$), which can also be phrased as follows:

\begin{definition}
\label{rrreinhardt} Reinhardt principle. \newline
There exists a non-trivial embedding $\jmath:{\rm V}\rightarrow{\rm V}$.
\end{definition}

This basic principle is inconsistent with ZFC. Kunen proved (in \cite{MR0311478}) that if $\jmath$ is an elementary embedding from ${\rm V}$ into ${\rm V}$, then $\jmath$ must be the identity. Several proofs are known (see \cite{MR1994835}, pp. 318-324), all of them employ the axiom of choice. The original proof is based on the following:

\begin{theorem}
\label{jjjj} $\omega$-J\'onsson algebras. \newline
For every infinite cardinal $\lambda$ there exists a function $f:[\lambda]^\omega\rightarrow\lambda$ such that $y\in[\lambda]^\lambda\Rightarrow f``[y]^\omega=\lambda$.
\end{theorem}

The proof of the combinatorial theorem appeared first in \cite{MR0209161}. In the common notation of the square brackets we may simply write $\lambda\nrightarrow[\lambda]^\omega_\lambda$. We indicate that the proof makes use of the axiom of choice, and no proof of the failure of Reinhardt principle in ZF is known.

But in the frame of ZFC, Reinhardt principle casts a limitation on the existence of too strong notions of large cardinals. In this paper we deal with two axioms, labeled as I1 and I2. From now on, an elementary embedding $\jmath$ is a non-trivial elementary embedding (i.e., $\jmath$ is not the identity):

\begin{definition}
\label{iiii} The axioms I1 and I2.
\begin{enumerate}
\item [$(\aleph)$] I1 means that for some ordinal $\delta$ there exists an elementary embedding $\jmath:{\rm V}_{\delta+1}\rightarrow{\rm V}_{\delta+1}$.
\item [$(\beth)$] I2 asserts that there is an elementary embedding $\jmath:{\rm V}\rightarrow M$ such that ${\rm V}_\delta\subseteq M$ for some $\delta>{\rm crit}(\jmath)$ satisfying $\jmath(\delta)=\delta$.
\end{enumerate}
\end{definition}

The axioms I1 and I2 were introduced, first, in \cite{MR0376347}. Interesting results which follow from these (and similar) axioms are proved in \cite{MR1954738}.
Magidor observed that if $\lambda\nrightarrow [\lambda]^{\aleph_0\text{-bd}}_\lambda$ for every (strong limit) $\lambda$, then I1 is refuted. Living in ${\rm V}_{\lambda+1}$, a function $f$ which exemplifies $\lambda\nrightarrow [\lambda]^{\aleph_0\text{-bd}}_\lambda$ plays the role of an $\omega$-J\'onsson function in the proof of Kunen, and the existence problem of such a function is attributed to Menachem Magidor (see \cite{MR1994835} Question 24.1). By a paper of Shioya (see \cite{MR1271428}) it refutes the axiom I2 as well.

Being the central notion of this paper, an appropriate name is in order:

\begin{definition}
\label{mmm} Magidor cardinals. \newline
Let $\lambda$ be an infinite cardinal.
\begin{enumerate}
\item [$(\aleph)$] A function $f:[\lambda]^{\aleph_0\text{-bd}}\rightarrow\alpha$ is a Magidor function for $\lambda$ (with $\alpha$-many colors) iff $f``[A]^{\aleph_0\text{-bd}}=\alpha$ whenever $A\in[\lambda]^\lambda$.
\item [$(\beth)$] A cardinal $\lambda$ with countable cofinality is a Magidor cardinal iff $\lambda\rightarrow[\lambda]^{\aleph_0\text{-bd}}_\lambda$.
\end{enumerate}
\end{definition}

If $\cf(\lambda)>\aleph_0$ then a Magidor function on $\lambda$ is simply an $\omega$-J\'onsson function, hence non-existent in ZFC. We shall prove that Magidor cardinals are large cardinals (in the philosophical sense, i.e. their existence is axiomatic and cannot be derived from ZFC). Moreover, we shall see that these cardinals (if exist) are situated in a fairly high position among their friends in the table of large cardinals.

We try to follow the standard notation. If $\lambda$ is a cardinal then $[\lambda]^{\aleph_0}$ is the collection of subsets of $\lambda$ whose size is $\aleph_0$. By $[\lambda]^{\aleph_0\text{-bd}}$ we denote the bounded subsets of $\lambda$ whose size is $\aleph_0$. Arrows notation with bd as a supercript are to be interpreted in the same manner. For example, $\lambda\nrightarrow [\lambda]^{\aleph_0\text{-bd}}_\alpha$ means that there exists a coloring $c:[\lambda]^{\aleph_0\text{-bd}}\rightarrow\alpha$ such that for every $y\in[\lambda]^\lambda$ we have $c``[y]^{\aleph_0\text{-bd}}=\alpha$. By the notation $\lambda\nrightarrow [\lambda]^{\aleph_0\text{-bd}}_{\theta,<\theta}$ we mean that the number of colors obtained is less than $\theta$.

For $\tau=\cf(\tau)<\lambda$ let $S^\lambda_\tau$ be the set $\{\delta<\lambda:\cf(\delta)=\tau\}$. Notice that if $\lambda$ is a regular cardinal then $S^\lambda_\tau$ is a stationary subset of $\lambda$. We shall use this notation even if $\lambda>\cf(\lambda)=\aleph_0$. A set $B$ is $\theta$-closed if the following requirement holds: for every increasing sequence $\langle\delta_i:i<\theta\rangle$ of members of $B$ we have $\delta={\rm sup}\{\delta_i:i<\theta\}\in B$.
If $\jmath$ is an elementary embedding between transitive models of ZFC then ${\rm crit}(\jmath)$ is the first ordinal moved by $\jmath$ and ${\rm fix}(\jmath)$ is the first fixed point of $\jmath$ above the critical point.

A word about square principles. The classical square has been defined by Jensen in his profound analysis of the constructible universe \cite{MR0309729}. Actually, Jensen phrased several versions of the square, and also described a straightforward way to force their existence. A recurrent question in this area is the relationship between squares and large cardinals. By and large, the square is an un-compact creature while large enough cardinals are compact. Consequently, the higher part of the large cardinals table and the square are orthogonal, but one can weaken the amount of square in order to incorporate larger cardinals in the universe, and this is the main theme in the current paper.
Recall that for a set $C$ we define ${\rm acc}(C)=\{\alpha\in C:\alpha={\rm sup}\{\beta:\beta\in C\cap\alpha\} \}$.
We shall force the existence of the following principle:

\begin{definition}
\label{ggglobalsq} Partial global square. \newline
Suppose $\delta=\cf(\delta)$. \newline
The $\delta$-partial global square means that there exists $S\supseteq S^{On}_\delta$ and a sequence $\langle C_\alpha:\alpha\in S,\alpha>\delta\rangle$ such that ${\rm otp}(C_\alpha)\leq\delta$ for every $\alpha$, ${\rm acc}(C_\alpha)\subseteq S$ for every $\alpha$ and $C_\alpha\cap\beta=C_\beta$ whenever $\beta\in{\rm acc}(C_\alpha)$.
\end{definition}

The historical background of Magidor cardinals goes back more than forty years ago, shortly after the concept of Reinhardt cardinals was introduced and Kunen's inconsistency has been proved. Due to Magidor himself, \cite{pc}, it emerged out of an innocent attempt to prove the inconsistency of I1 along the line of Kunen's proof. Things have been changed, notwithstanding. In the wake of Woodin's work about ${\rm AD}^{L(\mathbb{R})}$ and the axiom I0 (see \cite{MR2768702} and \cite{MR1994835} pp. 328-329) it seems that set theorists feel that the axioms of rank-into-rank are stable enough. In particular, the existence of Magidor cardinals is confident at least like the existence of finite ordinals.

We follow the notation of \cite{MR1940513} in general, with the following important exception. We use the Jerusalem forcing notation, i.e. $p\leq q$ means that $q$ is stronger than $p$. We shall use the forcing of Laver from \cite{MR0472529} by assuming that a supercompact cardinal $\kappa$ is indestructible under $\kappa$-directed-closed forcing notions.
For an excellent background regarding the contents of this paper we suggest \cite{MR1994835} and \cite{MR2768691}.

We are grateful to Menachem Magidor for a very helpful discussion on the subject of Magidor cardinals. We also thank the referee for a wonderful work which was much deeper than just proofreading. In particular, the current version of Theorem \ref{omegaclosed}, Proposition \ref{soduko} and Proposition \ref{solreferee} are due to the referee.

\newpage

\section{Combinatorial properties of Magidor cardinals}

We commence with a coding lemma, needed for several claims below:

\begin{lemma}
\label{ccodding} Coding reals.
\begin{enumerate}
\item [$(\aleph)$] There exists a function $r:[\omega]^\omega\rightarrow [\omega]^\omega$ such that for every $x,z\in[\omega]^\omega$ there exists a subset $y\in[x]^\omega$ for which $r(y)=z$.
\item [$(\beth)$] There exists a function $r:[\omega]^\omega\rightarrow 2^{\aleph_0}$ so that for every $x\in[\omega]^\omega$ and every $\beta<2^{\aleph_0}$ there exists a subset $y\in[x]^\omega$ such that $r(y)=\beta$.
\item [$(\gimel)$] Both assertions can be implemented on $[\lambda]^\omega$ in lieu of $[\omega]^\omega$ for every infinite cardinal $\lambda$, i.e. there is a function $r:[\lambda]^\omega\rightarrow 2^{\aleph_0}$ so that for every $x\in[\lambda]^\omega$ and every $\beta<2^{\aleph_0}$ there exists a subset $y\in[x]^\omega$ such that $r(y)=\beta$.
\end{enumerate}
In all cases we call $r$ a coding reals function.
\end{lemma}

\par\noindent\emph{Proof}. \newline
We prove the first part of the lemma, the proof of the second and third part is just the same. We fix an enumeration $\{x_\alpha:\alpha<2^{\aleph_0}\}$ of the members of $[\omega]^\omega$, in which every member appears $2^{\aleph_0}$-many times. Likewise, we fix an enumeration $\{z_\beta:\beta<2^{\aleph_0}\}$ of the members of $[\omega]^\omega$, in which every member appears only once. The function $r$ is defined by induction on $\alpha<2^{\aleph_0}$.

Arriving at $\alpha$, let $\beta$ be the order type of the set $\{\gamma<\alpha:x_\gamma=x_\alpha\}$. Choose a subset $y_\alpha\in[x_\alpha]^\omega$ such that $y_\alpha\neq y_\gamma$ for every $\gamma<\alpha$. The choice is possible since $\alpha<2^{\aleph_0}$ and we have $2^{\aleph_0}$ members of $[x_\alpha]^\omega$ at our disposal. Define $r(y_\alpha)=z_\beta$.

Assume now that $x,z\in[\omega]^\omega$. By the nature of our enumerations, there exists a unique ordinal $\beta$ so that $z=z_\beta$. Since $\beta<2^{\aleph_0}$ there exists an ordinal $\alpha<2^{\aleph_0}$ such that $x=x_\alpha$ and $\beta={\rm otp}(\{\gamma<\alpha:x_\gamma=x_\alpha\})$. By the $\alpha$-th stage of the construction there is some $y_\alpha\in[x]^\omega$ such that $r(y_\alpha)=z_\beta=z$, so we are done.

For getting the same coding with respect to $[\lambda]^\omega$, let $\mu=\lambda^\omega$ and enumerate the members of $[\lambda]^\omega$ in such a way that every set appears $2^{\aleph_0}$ many times. By the same process as above, one can define now the function $r$ which codes every real number.

\hfill \qedref{ccodding}

Our main goal is to refute the existence of Magidor cardinals from a weak version of the square principle. We are trying to employ the weakest square from the large cardinals point of view, i.e. a principle which can live happily with strong axioms of large cardinals. As a first step we show that one can replace the demand of $\lambda$-many colors (in the definition of Magidor cardinals) by $\alpha$-many colors for some $\alpha<\lambda$:

\begin{lemma}
\label{ccccolors} If $\lambda\nrightarrow[\lambda]^{\aleph_0\text{-bd}}_\alpha$ for every $\alpha<\lambda$ then $\lambda\nrightarrow[\lambda]^{\aleph_0\text{-bd}}_\lambda$.
\end{lemma}

\par\noindent\emph{Proof}. \newline
For every $\alpha<\lambda$ we choose a coloring $c_\alpha:[\lambda]^{\aleph_0\text{-bd}}\rightarrow\alpha$ which exemplifies the relation $\lambda\nrightarrow[\lambda]^{\aleph_0\text{-bd}}_\alpha$. Now we define a single coloring $c:[\lambda]^{\aleph_0\text{-bd}}\rightarrow\lambda$ as follows. If $x\in[\lambda]^{\aleph_0\text{-bd}}$ then let $\gamma$ be ${\rm min}(x)$, and set $c(x)=c_\gamma(x\setminus\{\gamma\})$.

Suppose $y\in[\lambda]^\lambda$ and $\beta<\lambda$. Choose an ordinal $\alpha\in y$ so that $\beta<\alpha$, and define $y'=y\setminus(\alpha+1)$. Clearly $y'\in[\lambda]^\lambda$, hence there exists some $x\in[y']^{\aleph_0\text{-bd}}$ such that $c_\alpha(x)=\beta$. Let $z$ be $x\cup\{\alpha\}$. It follows that $c(z)=c_\alpha(x)=\beta$. Since $z\in[y'\cup\{\alpha\}]^{\aleph_0\text{-bd}}\subseteq [y]^{\aleph_0\text{-bd}}$ we are done.

\hfill \qedref{ccccolors}

It follows from the above claim that if $\lambda$ is a Magidor cardinal then there exists a first ordinal $\alpha$ so that $\lambda\rightarrow[\lambda]^{\aleph_0\text{-bd}}_\alpha$, and we denote this ordinal by $\alpha_M$.
Similarly, if $\lambda$ is a J\'onsson cardinal then we let $\alpha_J$ be the first ordinal $\alpha$ so that $\lambda\rightarrow[\lambda]^{<\omega}_\alpha$.
Something further can be said about $\alpha_M$. The first part of the following claim is modelled after \cite{MR0337616}, and the second part applies to Magidor cardinals which emerge out of elementary embeddings. By $I1(\kappa,\lambda)$ we mean that $\jmath:{\rm V}_{\lambda+1}\rightarrow {\rm V}_{\lambda+1}$ is an elementary embedding so that $\kappa=\crit(\jmath)$, and hence $\lambda=fix(\jmath)$ is a Magidor cardinal.

\begin{claim}
\label{mmalpha} Some $\alpha_M$ properties. \newline
Let $\lambda$ be a Magidor cardinal.
\begin{enumerate}
\item [$(a)$] The ordinal $\alpha_M$ is a regular cardinal, and $\lambda\rightarrow[\lambda]^{\aleph_0\text{-bd}}_{\alpha_M,<\alpha_M}$.
\item [$(b)$] If $I1(\kappa,\lambda)$ then $\alpha_M<\kappa$.
\item [$(c)$] $\alpha_J\leq\alpha_M$.
\item [$(d)$] $2^{\aleph_0}<\alpha_M$.
\end{enumerate}
\end{claim}

\par\noindent\emph{Proof}. \newline
For part $(a)$ assume towards contradiction that $\cf(\alpha_M)<\alpha_M$. Fix any cofinal function $t:\cf(\alpha_M)\rightarrow\alpha_M$. Likewise, choose a function $f_\alpha:[\lambda]^{\aleph_0\text{-bd}}\rightarrow\alpha$ for every $\alpha<\alpha_M$ which exemplifies the negative relation $\lambda\nrightarrow[\lambda]^{\aleph_0\text{-bd}}_\alpha$. Denote the function $f_{\cf(\alpha_M)}$ by $g$.

Let $w$ be the set of infinite countable ordinals $\{\beta:\omega\leq\beta<\omega_1\}$. We choose a one-to-one mapping $p:w\times w\rightarrow w$ such that for every $(\beta_0,\beta_1)\in w\times w$ we have $\beta_0+\beta_1<p(\beta_0,\beta_1)$ (the sign $+$ refers to ordinals sum). Notice that $p$ need not be surjective (and usually, it is not).

Given a set $x\in[\lambda]^{\aleph_0\text{-bd}}$ let $\{x_j:j<\beta\}$ be an enumeration of the members of $x$ in increasing order. We define $h(x)$ as follows. If $\beta\notin{\rm rang}(p)$ then $h(x)=0$. If $\beta\in{\rm rang}(p)$ then there exists a unique pair $(\beta_0,\beta_1)$ so that $p(\beta_0,\beta_1)=\beta$. Since $\beta_0+\beta_1<\beta$ we can define safely the subsets $y=\{x_j:j<\beta_0\}$ and $z=\{x_j:\beta_0\leq j<\beta_1\}$. Define $h(x)=f_{t(g(y))}(z)$.

Assume $A\in[\lambda]^\lambda$ and $\gamma<\alpha_M$. First of all, we choose an ordinal $\alpha<\cf(\alpha_M)$ such that $\gamma<t(\alpha)<\alpha_M$. By the definition of $g$ there exists $y\in[A]^{\aleph_0\text{-bd}}$ such that $g(y)=\alpha$. Let $\beta_0$ be the order type of $(y,<)$.
By the definition of $f_{t(g(y))}$ there exists a subset $z\in[A-{\rm sup}(y)]^{\aleph_0\text{-bd}}$ such that $f_{t(g(y))}(z)=\gamma$. Let $\beta_1$ be the order type of $(z,<)$. Let $\beta$ be $p(\beta_0,\beta_1)$, and recall that $\beta_0+\beta_1<\beta$. We choose a subset $z^+=\{x_j:\beta_0+\beta+1\leq x_j<\beta\}\in[A-{\rm sup}(z)]^{\aleph_0\text{-bd}}$, and let $x=y\cup z\cup z^+$. It follows that $h(x)=f_{t(g(y))}(z)=\gamma$, hence $h$ exemplifies the negative relation $\lambda\nrightarrow[\lambda]^{\aleph_0\text{-bd}}_{\alpha_M}$, a contradiction.

The second part of the assertion (i.e. $\lambda\rightarrow[\lambda]^{\aleph_0\text{-bd}}_{\alpha_M,<\alpha_M}$) follows from the proof above. Indeed, for creating the contradictory $h$ we need only the fact that the range of $g$ is unbounded in $\alpha_M$, and $f_\alpha$ for each value of $g$.

For part $(b)$ assume that $\delta=\alpha_M\geq\kappa=\crit(\jmath)$. Since $\delta<\lambda=fix(\jmath)$ we know that $\jmath(\delta)>\delta$. Choose a function $f:[\lambda]^{\aleph_0\text{-bd}}\rightarrow\delta$ which exemplifies the Magidority of $\lambda$. By elementarity, $\jmath f:[\lambda]^{\aleph_0\text{-bd}}\rightarrow\jmath(\delta)$ is a Magidor function, and it lies in ${\rm V}_{\lambda+1}$. However, $\jmath(\delta)$ is $\alpha_M$ by elementarity, which is an absurd since $\jmath(\delta)>\delta$.

Moving to $(c)$, suppose we are given a Magidor cardinal $\lambda$ and fix an $\omega$-sequence of ordinals $\langle\alpha_n:n\in\omega\rangle$ which tends to $\lambda$. For every $x\in[\lambda]^{\aleph_0\text{-bd}}$ let $n_x$ be the first natural number for which $x\cap\alpha_{n_x}\neq\emptyset$. Define:

\begin{displaymath}
m(x)= \left\{
\begin{array}{ll} |x\cap\alpha_{n_x}| & \textrm{if the intersection is finite}\\
0 & \textrm{if the intersection is infinite}
\end{array} \right.
\end{displaymath}

Choose any ordinal $\beta<\alpha_J$. We shall construct a function $f:[\lambda]^{\aleph_0\text{-bd}}\rightarrow\beta$ which omits no color, hence proving that $\beta<\alpha_M$ for every $\beta<\alpha_J$. To begin with, let $g:[\lambda]^{<\omega}\rightarrow\beta$ be a J\'onsson function in the common sense. Suppose $x=\{x_n:n\in\omega\}\in[\lambda]^{\aleph_0\text{-bd}}$, and let $z$ be the finite subset $\{x_{m(x)},\ldots,x_{2m(x)-1}\}$. We let $f(x)=g(z)$.

For proving that $f$ is as required, assume $A\in[\lambda]^\lambda$ and $\gamma<\beta$ is any color. Let $n$ be the first natural number for which $|A\cap[\alpha_n,\alpha_{n+1})|\geq\aleph_0$.
By the J\'onssonicity of $g$ we can choose some $\eta=\eta_\gamma\in[A\setminus\alpha_{n+1}]^{<\omega}$ such that $g(\eta)=\gamma$.

Pick any subset of $A\cap[\alpha_n,\alpha_{n+1})$ of size $|\eta|$, say $\eta_0$. Choose a subset $\eta_1\in[A]^{\aleph_0\text{-bd}}$ so that ${\rm max}(\eta)<{\rm min}(\eta_1)$. Set $x=\eta_0\cup\eta\cup\eta_1$. By the definition of $f$ we have $f(x)=g(\eta)=\gamma$, so we are done.

Finally, fix a coding reals function $r_\alpha:[\alpha]^\omega\rightarrow 2^{\aleph_0}$ for every $\alpha\in[\omega,\lambda)$. Given $x\in[\lambda]^{\aleph_0\text{-bd}}$ such that ${\rm otp}(x)=\omega$, let $\alpha$ be ${\rm sup}(x)$. Define $f:[\lambda]^{\aleph_0\text{-bd}}\rightarrow 2^{\aleph_0}$ as follows. If ${\rm otp}(x)\neq\omega$ then $f(x)=0$ and if ${\rm otp}(x)=\omega$ then $f(x)=r_\alpha(x)$. By the coding reals lemma, $f$ exemplifies the negative relation $\lambda\nrightarrow [\lambda]^{\aleph_0\text{-bd}}_{2^{\aleph_0}}$ and even $\lambda\nrightarrow [\omega]^{\aleph_0\text{-bd}}_{2^{\aleph_0}}$, so we are done.

\hfill \qedref{mmalpha}

The third part of the above claim yields the following:

\begin{theorem}
\label{mt1} Magidor cardinals and J\'onsson cardinals. \newline
Suppose $\lambda$ is a Magidor cardinal. \newline
Then $\lambda$ is a J\'onsson cardinal as well. \newline
In particular, there are no Magidor cardinals in the constructible universe.
\end{theorem}

\par \noindent \emph{Proof}. \newline
Assume $\lambda$ is not a J\'onsson cardinal. For every $\beta<\lambda$ let $f_\beta$ exemplify the negative relation $\lambda\nrightarrow[\lambda]^{\aleph_0\text{-bd}}_\beta$. Such a function exists by virtue of part $(c)$ of Claim \ref{mmalpha}. However, it follows now from Lemma \ref{ccccolors} that $\lambda$ is not a Magidor cardinal, a contradiction.

\hfill \qedref{mt1}

\begin{remark}
\label{refereer} Part $(c)$ of Claim \ref{mmalpha}, used in the above proof, can be proved easily by the following argument. Suppose $\alpha<\alpha_J$ and fix a function $f:[\lambda]^{<\omega}\rightarrow\alpha$ which exemplifies $\lambda\nrightarrow[\lambda]^{<\omega}_\alpha$.

Define $g:[\lambda]^{\aleph_0\text{-bd}}\rightarrow\alpha$ as follows. If ${\rm otp}(x)=\omega+n$ and $\eta$ is the $\omega$-th member of $x$ then $g(x)=f(x\setminus\eta)$. Otherwise, $g(x)=0$. Clearly, $g$ exemplifies $\lambda\nrightarrow[\lambda]^{\aleph_0\text{-bd}}_\alpha$, so $\alpha<\alpha_M$.

This simple proof has been suggested by the referee of this paper, and we thank him or her. We keep, however, the indirect proof given above as it might be useful for other versions of Magidority which are dictated in Definition \ref{oordertypemagidor} below.
\end{remark}

\hfill \qedref{refereer}

The last part of Claim \ref{mmalpha} yields another interesting consequence:

\begin{corollary}
\label{bbbb} Magidority and the continuum. \newline
If $\lambda$ is a Magidor cardinal then $\lambda>2^{\aleph_0}$, while if there exists a J\'onsson cardinal then it is consistent that there is a J\'onsson cardinal below the continuum.
\end{corollary}

\par\noindent\emph{Proof}. \newline
The first part of the assertion follows from the fact that $2^{\aleph_0}<\alpha_M<\lambda$.
The second part appears in \cite{MR0363906}, and we describe shortly the argument for completeness. Assume that $\lambda$ is a J\'onsson cardinal, and let $\kappa=\lambda^+$. Let $\mathbb{P}$ be ${\rm Add}(\omega,\kappa)$, and force with $\mathbb{P}$ in order to add $\kappa$-many Cohen reals.

Clearly, ${\rm V}^{\mathbb{P}}\models 2^{\aleph_0}=\kappa$. We claim that $\lambda$ is a J\'onsson cardinal in the generic extension ${\rm V}^{\mathbb{P}}$. For this end, assume that $c:[\lambda]^{<\omega}\rightarrow\lambda$ is a coloring in ${\rm V}^{\mathbb{P}}$. Let $\name{c}$ be a name of this coloring in ${\rm V}$. The value of $\name{c}$ is determined by a small set of conditions, since $\mathbb{P}$ is $ccc$. It follows that the J\'onssonicity of $\lambda$ from the ground model is preserved, as we can define a function $f\in{\rm V}$ to exemplify it, so we are done.

\hfill \qedref{bbbb}

The above corollary is suggestive also for Magidor cardinals, in the sense that one can force them to be below $2^{\aleph_1}$. It means that although Magidor cardinals are large cardinals in the philosophical sense (i.e., their existence is axiomatic), if there exists a Magidor cardinal one can reduce its magnitude in the $\beth$-scale. In particular, a Magidor cardinal need not be a limit of strongly inaccessible cardinals, and can be smaller than the first strongly inaccessible:

\begin{claim}
\label{ssmallmagidor} Small Magidor cardinals. \newline
If there is a Magidor cardinal then it is consistent that there is a Magidor cardinal below $2^{\aleph_1}$.
\end{claim}

\par\noindent\emph{Proof}. \newline
We shall prove the following general assertion. Suppose $\mathbb{P}$ is $\aleph_1$-complete and $\alpha_M$-cc. If $G\subseteq\mathbb{P}$ is a generic subset, then $\lambda$ is still a Magidor cardinal in $V[G]$ and $\alpha_M^{{\rm V}[G]}\leq \alpha_M^{{\rm V}}$.
Let $\alpha_M$ be $\alpha_M^{{\rm V}}$ and assume that $\name{c}:[\lambda]^{\aleph_0\text{-bd}}\rightarrow\alpha_M$ is a name of a coloring. We have to find $\name{A}\in[\lambda]^\lambda$ for which $\Vdash_{\mathbb{P}} \name{c}``[\name{A}]^{\aleph_0\text{-bd}}\neq \check{\alpha}_M$.
For this end, we define in ${\rm V}$ a function $g:[\lambda]^{\aleph_0\text{-bd}}\rightarrow\alpha_M$, and by the Magidority of $\lambda$ in ${\rm V}$ we can choose a subset $A\in[\lambda]^\lambda$ for which $|g``[A]^{\aleph_0\text{-bd}}|<\alpha_M$. We shall see that $\name{c}$ omits colors on $\check{A}$.

Given any set $y\in[\lambda]^{\aleph_0\text{-bd}}$ we know that $\name{c}(y)$ is a name of an ordinal in $\alpha_M$. Since $\mathbb{P}$ is $\alpha_M$-cc there is an antichain $\mathcal{A}_y$ of size strictly less than $\alpha_M$ which forces a value to $\name{c}(y)$. Define:

$$
\delta(y)={\rm sup}\{\delta:\exists p\in\mathcal{A}_y, p \Vdash \name{c}(y)=\delta\}.
$$

Clearly, $\delta(y)<\alpha_M$. Set $g(y)=\delta(y)$ for every $y\in[\lambda]^{\aleph_0\text{-bd}}$.
Since the forcing relation is definable in ${\rm V}$ we have $g\in{\rm V}$, and $g:[\lambda]^{\aleph_0\text{-bd}}\rightarrow\alpha_M$. Choose $A\in[\lambda]^\lambda$ so that $|g``[A]^{\aleph_0\text{-bd}}|<\alpha_M$, and let $\beta$ be ${\rm sup}(g``[A]^{\aleph_0\text{-bd}})$. Now $\check{A}\in {\rm V}[G]$, and if $\name{y}\in[\check{A}]^{\aleph_0\text{-bd}}$ then $y\in[A]^{\aleph_0\text{-bd}}$ by the $\aleph_1$-completeness of $\mathbb{P}$, so $\Vdash_{\mathbb{P}} \name{c}(y)\leq g(y)\leq\beta$. It follows that $\lambda$ is a Magidor cardinal in ${\rm V}[G]$, as required.

Let $\lambda$ be Magidor in ${\rm V}$, and apply the general assertion at the beginning of the proof to the Cohen forcing for adding $\lambda^+$-many subsets of $\aleph_1$. This forcing satisfies the above requirement from $\mathbb{P}$, (notice that $\mathbb{P}$ is $(2^{\aleph_0})^+$-cc so also $\alpha_M$-cc since $2^{\aleph_0}<\alpha_M$) hence $\lambda$ is still a Magidor cardinal in the forcing extension. However, $\lambda<2^{\aleph_1}$, so we are done.

\hfill \qedref{ssmallmagidor}

The first part of Claim \ref{mmalpha} shows that a Magidor cardinal has some Rowbottom properties when the number of colors is $\alpha_M$. We can prove a parallel result for every $\alpha\in[\alpha_M,\lambda)$ provided that $\alpha$ has uncountable cofinality:

\begin{theorem}
\label{omegaclosed} Assume $\alpha$ is a cardinal and $\cf(\alpha)>\omega$. Let $\lambda$ be a Magidor cardinal, and assume that $\lambda\rightarrow [\lambda]^{\aleph_0\text{-bd}}_{\alpha^{\aleph_0}}$.
\begin{enumerate}
\item [$(a)$] $\lambda\rightarrow [\lambda]^{\aleph_0\text{-bd}}_{\alpha^{\aleph_0},<\alpha}$.
\item [$(b)$] If $\lambda\rightarrow [\lambda]^{\aleph_0\text{-bd}}_\alpha$ then $\lambda\rightarrow [\lambda]^{\aleph_0\text{-bd}}_{\alpha,<\alpha}$.
\item [$(c)$] If $\alpha<\alpha_M$ then $\alpha^{\aleph_0}<\alpha_M$.
\end{enumerate}
\end{theorem}

\par\noindent\emph{Proof}. \newline
Let $f$ be a function from $[\lambda]^{\aleph_0\text{-bd}}$ into $\alpha^{\aleph_0}$. We have to find a set $A\in[\lambda]^\lambda$ such that $|f``[A]^{\aleph_0\text{-bd}}|<\alpha$. For this end, fix an $\omega$-J\'onsson function $g:[\alpha^{\aleph_0}]^\omega\rightarrow\alpha^{\aleph_0}$. Fix also an enumeration $\langle b_\zeta:\zeta<\alpha^{\aleph_0}\rangle$ of all the members of $[\alpha^{\aleph_0}]^\omega$ without repetitions.

Notice that $|{}^{\omega\cdot\omega}([\omega_1]^\omega)|=2^{\aleph_0}$. By the coding reals lemma we can choose for any $\delta<\lambda$ a map $r_\delta:[\delta]^\omega\rightarrow {}^{\omega\cdot\omega}([\omega_1]^\omega)$ such that for every $x\in[\delta]^\omega$ and every $s\in {}^{\omega\cdot\omega}([\omega_1]^\omega)$ there is $y\in[x]^\omega$ for which $r_\delta(y)=s$.

We define another function $h$ from $[\lambda]^{\aleph_0\text{-bd}}$ into $\alpha^{\aleph_0}$. Given any $x\in[\lambda]^{\aleph_0\text{-bd}}$ we ask whether there exists a limit ordinal $\beta<\omega_1$ such that ${\rm otp}(x)=\beta+\omega$. If the answer is negative then $h(x)=0$. If the answer is positive, we decompose $x$ in the following way. First, let $\langle x_i:i<\beta+\omega\rangle$ be an increasing enumeration of the members of $x$. Denote ${\rm sup}(x)$ by $\delta$ and set $s=r_\delta(\{x_{\beta+i}:i<\omega\})$. We may assume that $s=\langle a_j:j<\omega\cdot\omega\rangle$ when $a_j\in[\beta]^\omega$ for every $j<\omega\cdot\omega$.

For every $n,j\in\omega$ let $y^n_j=\{x_i:i\in a_{\omega\cdot n+j}\}$. For every $n\in\omega$ let $z_n=\{f(y^n_j):j\in\omega\}$. By the definition of $f$ we have $z_n\in[\alpha^{\aleph_0}]^\omega$. Hence there exists a unique ordinal $\zeta_n<\alpha^{\aleph_0}$ such that $z_n=b_{\zeta_n}$. We define $h(x)=g(\{\zeta_n:n\in\omega\})$. This accomplishes the definition of $h$.

We choose a set $A\in[\lambda]^\lambda$ such that $h``[A]^{\aleph_0\text{-bd}}\neq\alpha^{\aleph_0}$. We claim that $|f``[A]^{\aleph_0\text{-bd}}|<\alpha$. For proving this claim, assume toward contradiction that $|f``[A]^{\aleph_0\text{-bd}}|\geq\alpha$, and let $\eta<\alpha^{\aleph_0}$ be any ordinal. We shall find some $x\in[A]^{\aleph_0\text{-bd}}$ for which $h(x)=\eta$, thus proving that $h``[A]^{\aleph_0\text{-bd}}=\alpha^{\aleph_0}$, a contradiction.

Since $\cf(\alpha)>\omega$ and $\cf(\lambda)=\omega$ we can find $\mu<\lambda$ such that $|f``[A\cap\mu]^{\aleph_0\text{-bd}}|\geq\alpha$. Let $B=f``[A\cap\mu]^{\aleph_0\text{-bd}}$ and $C=\{\zeta<2^{\aleph_0}:b_\zeta\in[B]^\omega\}$. The cardinality of $C$ is $\alpha^{\aleph_0}$ as $|B|\geq\alpha$, and hence $g``[C]^\omega=\aleph^{\aleph_0}$. In particular, we can choose an element $z\in[C]^\omega$ so that $g(z)=\eta$.

Let $\langle\zeta_n:n\in\omega\rangle$ be an increasing enumeration of the members of $z$. By the definition of $C$, for every $n\in\omega$ we have $b_{\zeta_n}\in[B]^\omega$. For every $n,j\in\omega$ we find $y^n_j\in[A\cap\mu]^\omega$ such that $b_{\zeta_n}=\{f(y^n_j):j\in\omega\}$. Choose any $y\in[A]^{\aleph_0\text{-bd}}$ such that $\mu<{\rm min}(y)$ and ${\rm otp}(\bigcup\limits_{n,j\in\omega}y^n_j\cup y)=\beta$ for some limit ordinal $\beta<\omega_1$.

Let $x'$ be $\bigcup\limits_{n,j\in\omega}y^n_j\cup y$ and enumerate the members of $x'$ in increasing order by $\langle x_i:i<\beta\rangle$. As above, let $\langle a_j:j<\omega\cdot\omega\rangle$ be the corresponding decomposition, so that $a_j\in[\beta]^\omega$ and $y^n_j=\{x_i:i\in a_{\omega\cdot n+j}\}$. In order to get the correct order type we choose a set $w\in[A]^{\aleph_0\text{-bd}}$ such that ${\rm sup}(x')<{\rm min}(w), {\rm otp}(w)=\omega$, and $r_{{\rm sup}(w)}(w)=\langle a_j:j<\omega\cdot\omega\rangle$. Define $x=x'\cup w$, so $x\in[A]^{\aleph_0\text{-bd}}$ and notice that $h(x)=\eta$, so we are done.

\hfill \qedref{omegaclosed}

\begin{remark}
\label{rrr} It follows from the above theorem that every Magidor cardinal $\lambda$ is $\omega$-closed, i.e. $\alpha<\lambda\Rightarrow \alpha^{\aleph_0}<\lambda$. A slightly different proof of this fact appears in \cite{MR3666820}.
\end{remark}

\hfill \qedref{rrr}

Based on the above remark, one can show the following:

\begin{proposition}
\label{soduko} If there exists a Magidor cardinal $\lambda$ then there is a generic extension in which $\lambda$ is still Magidor and $\alpha_M=\aleph_2$.
\end{proposition}

\par\noindent\emph{Proof}. \newline
First observe that if $\lambda$ is Magidor, $\alpha=\cf(\alpha)<\lambda$ and $\lambda\rightarrow[\lambda]^{\aleph_0\text{-bd}}_\alpha$ then for every forcing notion $\mathbb{P}$ which is $\aleph_1$-complete and $\alpha$-cc we have $\Vdash_{\mathbb{P}}\lambda\rightarrow[\lambda]^{\aleph_0\text{-bd}}_\alpha$.

For proving this statement recall that if $\lambda\rightarrow[\lambda]^{\aleph_0\text{-bd}}_\alpha$ and $\alpha$ is regular then $\lambda\rightarrow[\lambda]^{\aleph_0\text{-bd}}_{\alpha,<\alpha}$. Now fix a condition $p\in\mathbb{P}$ and let $\name{f}$ be a name of a function from $[\lambda]^{\aleph_0\text{-bd}}$ into $\alpha$. For every $x\in[\lambda]^{\aleph_0\text{-bd}}$ let $R_x=\{\beta<\alpha:\exists q\geq p, q\Vdash\name{f}(x)=\check{\beta}\}$. Since $\mathbb{P}$ is $\alpha$-cc, $|R_x|<\alpha$. Define $g(x)={\rm sup}(R_x)$ and notice that $g(x)<\alpha$ since $\alpha$ is regular.

Since the forcing relation is definable in ${\rm V}, g\in{\rm V}$. By the construction, $p\Vdash\name{f}(x)\leq g(x)$ whenever $x\in[\lambda]^{\aleph_0\text{-bd}}$. Choose $A\in[\lambda]^\lambda$ so that $|g``[A]^{\aleph_0\text{-bd}}|<\alpha$. It follows that $p\Vdash {\rm sup}(\name{f}``[A]^{\aleph_0\text{-bd}})\leq {\rm sup}(g``[A]^{\aleph_0\text{-bd}})<\alpha$ (here we use the $\aleph_1$-completeness of $\mathbb{P}$) and hence $p\Vdash\name{f}``[A]^{\aleph_0\text{-bd}}\neq\alpha$.

Given a Magidor cardinal $\lambda$, let $\alpha$ be $((\alpha_M)^{\aleph_0})^+$. Since $\lambda$ is $\omega$-closed and limit, $\alpha<\lambda$. Likewise, $\alpha$ itself is $\omega$-closed and regular. Consequently, there is a forcing notion $\mathbb{P}$ which is $\aleph_1$-complete and $\alpha$-cc such that $\Vdash_{\mathbb{P}}\alpha=\aleph_2$. By the above argument, $\lambda\rightarrow[\lambda]^{\aleph_0\text{-bd}}_\alpha$ in the generic extension, and hence $\lambda$ is Magidor with $\alpha_M\leq\alpha$ in ${\rm V}^{\mathbb{P}}$. On the other hand, $\alpha_M\geq\alpha=\aleph_2$ since $\alpha_M>2^{\aleph_0}$ and hence $\Vdash_{\mathbb{P}}\alpha_M=\aleph_2$ as desired.

\hfill \qedref{soduko}

Magidor cardinals were defined with respect to $\omega$-bounded sets, regardless of the order type of these sets. The following definition is more sensitive:

\begin{definition}
\label{oordertypemagidor} $\beta$-Magidority. \newline
Let $\beta$ be an infinite ordinal.
\begin{enumerate}
\item [$(\aleph)$] $[\lambda]^{<\beta\text{-bd}}$ is the collection of bounded subsets of $\lambda$ whose order type is strictly less than $\beta$.
\item [$(\beth)$] $\lambda\rightarrow[\lambda]^{<\beta\text{-bd}}_\lambda$ iff for every $c:[\lambda]^{<\beta\text{-bd}}\rightarrow\lambda$ there exists $A\in[\lambda]^\lambda$ for which $c\upharpoonright[A]^{<\beta\text{-bd}}\neq\lambda$.
\item [$(\gimel)$] A cardinal $\lambda$ is $\beta$-Magidor iff $\lambda\rightarrow[\lambda]^{<\beta\text{-bd}}_\lambda$.
\item [$(\daleth)$] We call $\lambda$ strongly-Magidor iff $\lambda$ is $\beta$-Magidor for every $\beta<\lambda$.
\end{enumerate}
\end{definition}

Listed below are some basic observations, the proof of which is similar to the above proofs for the common Magidority (i.e., the case of $\beta=\omega_1$). We denote by $\alpha_M(\beta)$ the first ordinal $\alpha$ so that $\lambda\rightarrow[\lambda]^{<\beta\text{-bd}}_\alpha$.
\begin{itemize}
\item $2^{\aleph_0}<\alpha_M(\omega+1)$.
\item $\beta<\gamma\Rightarrow\alpha_M(\beta)\leq\alpha_M(\gamma)$.
\item If $I1(\kappa,\lambda)$ then $\lambda$ is strongly Magidor.
\end{itemize}

The proof of the last item begins like the proof of Magidority for I1 cardinals.
First we show that $\lambda$ is $<\theta$-Magidor for every $\theta<\kappa$.
Indeed, if $I1(\kappa,\lambda)$ and $f:[\lambda]^{<\theta\text{-bd}}\rightarrow\lambda$ exemplifies $\lambda\nrightarrow[\lambda]^{<\theta\text{-bd}}_\lambda$, then $\jmath f:[\lambda]^{<\theta\text{-bd}}\rightarrow\lambda$ exemplifies $\lambda\nrightarrow[\lambda]^{<\theta\text{-bd}}_\lambda$ by elementarity (we use the fact that $f$ belongs to ${\rm V}_{\lambda+1}$, as the sets in the domain of $f$ are bounded). However, $\jmath``\lambda\in[\lambda]^\lambda$ and $\kappa\notin\jmath f``[\jmath``\lambda]^{<\theta\text{-bd}}$.

Second, we can define below $\lambda$ elementary embeddings for which the critical point is more and more large, up to $\lambda$. For every $n\in\omega$ there is $\imath_n:{\rm V}_{\lambda+1}\rightarrow {\rm V}_{\lambda+1}$ such that $\crit(\imath_n)=\jmath^n(\kappa)$. This can be done by defining $\imath_n$ over ${\rm V}_\lambda$ and then extending it (essentially, in a unique way) to ${\rm V}_{\lambda+1}$ as described in \cite{MR1994835}, pp. 325-326. Now if we choose a sequence $\langle\theta_n:n\in\omega\rangle$ cofinal in $\lambda$ then we can use the first step described above in order to show that $\lambda$ is $<\theta_n$-Magidor for every $n\in\omega$.

\newpage

\section{Partial squares and Magidor cardinals}

Diminishing the number of colors in Lemma \ref{ccccolors}, we can phrase a useful combinatorial property which enables us to prove negative square brackets relations with respect to $\omega$-bounded subsets. This property is a special kind of reflection for certain stationary sets:

\begin{definition}
\label{ssss} Quilshon (pitchfork). \newline
Assume $\lambda>\delta=\cf(\delta)$. \newline
We say that $\pitchfork(\lambda,\delta)$ (or $\pitchfork_{\lambda,\delta}$) holds iff there is a collection $\{S_\gamma:\gamma<\delta\}$ of disjoint subsets of $\lambda$ so that $S_\gamma\cap\eta$ is a stationary subset of $\eta$ for every ordinal $\eta<\lambda$ with $\cf(\eta)=\delta$ and every $\gamma<\delta$.
\newline
We may replace $\lambda$ by any unbounded subset $S\subseteq\lambda$, so $\pitchfork_{S,\delta}$ means that we decompose $S$ rather than $\lambda$.
\end{definition}

The following theorem draws a connection between $\pitchfork(\lambda,\delta)$ and $\omega$-bounded J\'onsson functions:

\begin{theorem}
\label{mt2} Non-Magidority and the quilshon principle. \newline
Assume that $\lambda>\cf(\lambda)=\aleph_0$ is a Magidor cardinal, $S=S^\lambda_\omega$ and let $\alpha=\alpha_M<\lambda$ be the first ordinal so that $\lambda\rightarrow[\lambda]^{\aleph_0\text{-bd}}_\alpha$. \newline
Then $\neg(\pitchfork_{S,\delta})$ for every $\delta\in {\rm Reg}\cap[\alpha,\lambda)$.
\end{theorem}

\par\noindent\emph{Proof}. \newline
Assume to the contrary that $\pitchfork_{S,\delta}$ for some $\delta\in {\rm Reg}\cap[\alpha,\lambda)$. We shall prove that in this case $\lambda\nrightarrow[\lambda]^{\aleph_0\text{-bd}}_\delta$. For this end, let $\{S_\gamma:\gamma<\delta\}$ exemplify $\pitchfork_{S,\delta}$. For each ordinal $\eta<\lambda$ with $\cf(\eta)=\delta$ let $S_{\eta,\gamma}$ be $S_\gamma\cap\eta$. By the definition of the pithcfork, $\{S_{\eta,\gamma}:\gamma<\delta\}$ forms a partition of $S^\eta_\omega$ into $\delta$-many disjoint stationary sets.
Given $x\in[\lambda]^{\aleph_0\text{-bd}}$ let $\eta_x$ be the first ordinal such that $\cf(\eta_x)=\delta$ and $\eta_x>{\rm sup}(x)$. Since the cofinality of ${\rm sup}(x)$ is $\omega$, there exists a unique ordinal $\gamma\in\delta$ such that ${\rm sup}(x)\in S_{\eta_x,\gamma}$. Define $f(x)$ to be this ordinal.

Suppose $A\in[\lambda]^\lambda$ and let $\{a_\beta:\beta<\lambda\}$ be an increasing enumeration of the members of $A$. Denote the initial segment $\{a_\beta:\beta<\delta\}$ by $A_\delta$ and ${\rm sup}(A_\delta)$ by $\eta$, so $\cf(\eta)=\delta$. Let $\varepsilon<\delta$ be any color, let $A_\delta^{\omega-c\ell}$ be the $\omega$-closure of $A_\delta$, and notice that it meets every stationary subset of $S^\eta_\omega$. In particular, $S_{\eta,\varepsilon}\cap A_\delta^{\omega-c\ell}\neq\emptyset$. This means that for some $x\in[A_\delta]^\omega$ we have ${\rm sup}(x)\in S_{\eta,\varepsilon}$.

However, $A_\delta$ is bounded in $\lambda$ by the ordinal $\eta$, so $x\in[A]^{\aleph_0\text{-bd}}$. Likewise, if $\eta_x\leq\eta$ is the first ordinal above ${\rm sup}(x)$ with cofinality $\delta$ then ${\rm sup}(x)\in S_{\eta_x,\varepsilon}$ by the quilshon. It follows from the definition of $f$ and the fact that ${\rm sup}(x)\in S_{\eta_x,\varepsilon}= S_\varepsilon\cap\eta_x = (S_\varepsilon\cap\eta)\cap\eta_x= S_{\eta,\varepsilon}\cap\eta_x$ that $f(x)=\varepsilon$. Since $\varepsilon<\alpha$ was arbitrary we infer that $\lambda\nrightarrow[\lambda]^{\aleph_0\text{-bd}}_\delta$, a contradiction.

\hfill \qedref{mt2}

The next stage is essentially to prove that the quilshon follows from a partial global square, and even less:

\begin{claim}
\label{qqqilshonfromsaure} The quilshon claim. \newline
Assume $\delta=\cf(\delta)>\aleph_0$ and the $\delta$-partial global square principle holds. \newline
Then $\pitchfork_{\lambda,\delta}$ holds for every $\lambda>\delta$.
\end{claim}

\par\noindent\emph{Proof}. \newline
Fix any $\delta$-partial global square sequence of the form $\langle C_\alpha:\alpha\in S,\alpha>\delta\rangle$, when $S\supseteq S^{On}_\delta$.
We shall define a function $f:S\rightarrow\delta+1$, and then create a quilshon out of $f$.

To begin with, choose any partition of $\delta$ into $\delta$-many (disjoint) stationary sets $\langle T_\gamma:\gamma<\delta\rangle$. For every $\alpha<\delta$ let $f(\alpha)=\gamma$ iff $\alpha\in T_\gamma$. Notice that $T_\gamma=f^{-1}(\{\gamma\})$ for every $\gamma<\delta$, and a similar property will be maintained along the inductive construction of $f$. We would like to define $f$ also on $S\setminus\delta$, so assume $\alpha\in S\setminus\delta$ and distinguish two cases. If $\cf(\alpha)\geq\delta$ then let $f(\alpha)=\delta$, in which case $\alpha$ is uninteresting from our point of view. If $\cf(\alpha)<\delta$ then $f(\alpha)$ is defined as $f({\rm otp}(C_\alpha))$.

Ahead of defining the quilshon, we should prove that $f$ is well defined. The problem is non-existent when $\cf(\alpha)\geq\delta$, and if ${\rm otp}(C_\alpha)<\delta$ then $f({\rm otp}(C_\alpha))$ is well-defined by the initial decomposition into stationary sets. By the definition of partial global square, all the cases are covered.

Let $\lambda$ be an infinite cardinal above $\delta$ (typically, $\lambda$ is a singular cardinal of coubtable cofinality). Let $S_\gamma=\{\beta<\lambda:f(\beta)=\gamma\}$ for every $\gamma<\delta$. We claim that the collection $\langle S_\gamma:\gamma<\delta\rangle$ forms a quilshon for $\lambda$. For showing this, assume $\eta<\lambda$ and $\cf(\eta)=\delta$. Assume further that $\gamma<\delta$. We need to show that $S_\gamma\cap\eta$ is a stationary subset of $\eta$.

Fix any club $D_\eta\subseteq\eta$. We shall prove that $D_\eta\cap S_\gamma\neq \emptyset$. This can be done simply by showing that there exists an ordinal $\beta\in D_\eta$ such that $f(\beta)=\gamma$. We may assume that $D_\eta\subseteq{\rm acc}(C_\eta)$ (recall that $\delta$ is uncountable). Denote $\zeta={\rm otp}(C_\eta)$. By the global square properties, $\zeta<\eta$. Inasmuch as $\cf(\zeta)=\cf(\eta)=\delta$ we can use the induction hypothesis to conclude that $S_\gamma\cap\zeta$ is a stationary subset of $\zeta$.

Let $t:\zeta\rightarrow C_\eta$ be an increasing enumeration of the members of $C_\eta$. We translate the club $D_\eta$ into $D_\zeta=\{\gamma<\zeta:t(\gamma)\in D_\eta\}$, so $D_\zeta$ is a club in $\zeta$. By virtue of the induction hypothesis we can choose an ordinal $\zeta'\in D_\zeta\cap S_\gamma$. It means that $f(\zeta')=\gamma$, and $\beta=t(\zeta')\in D_\eta$. However, $f(\beta)=f({\rm otp}(C_\beta))$ which by coherence equals $f({\rm otp}(C_\eta\cap\beta))=f(\zeta')=\gamma$, so we are done.

\hfill \qedref{qqqilshonfromsaure}

Our main purpose is to prove the consistency of large cardinals with the fact that no Magidor cardinal exists. We begin with a limit of measurable cardinals. Recall that if $\lambda>\cf(\lambda)$ is a limit of measurable cardinals then $\lambda$ is J\'onsson (by results of Prikry, see \cite{MR0262075}). Nevertheless, it is consistent that there are many measurable cardinals while there are no Magidor cardinals at all.

Before proving it, recall that $\kappa$ is $1$-extendible iff there exists a cardinal $\lambda>\kappa$ and an elementary embedding $\jmath:\mathcal{H}(\kappa^+)\rightarrow\mathcal{H}(\lambda^+)$ so that $\crit(\jmath)=\kappa$ and $\jmath(\kappa)=\lambda$ (this definition comes from \cite{MR1942302}, and other equivalent formulations exist in the literature). It is easy to see that if $\kappa$ is $1$-extendible then there is a normal ultrafilter $U$ on $\kappa$ which concentrates on the measurable cardinals below $\kappa$ (see, e.g. \cite{MR1994835}, Proposition 23.1).

It is known that the existence of a global square implies the partial global square at every regular $\delta$.
Theorem 6.5 of \cite{MR1942302} asserts that if there exists some $1$-extendible cardinal $\kappa$ in the ground model (and an inaccessible cardinal above $\lambda$ where $\lambda$ is the target of the embedding) then one can force a global square while preserving the fact that $\kappa$ is $1$-extendible:

\begin{corollary}
\label{ccumandshimmer} Non-Magidority and limit of measurables. \newline
It is consistent that $\kappa$ is $1$-extendible and there are no Magidor cardinals. \newline
Consequently, it is consistent that there is a class of measurable cardinals and no Magidor cardinals.
In particular, a limit of measurable cardinals need not be a Magidor cardinal.
\end{corollary}

\par\noindent\emph{Proof}. \newline
Combine the above mentioned theorem from \cite{MR1942302} with Claim \ref{qqqilshonfromsaure} in order to get a model with $1$-extendible cardinal $\kappa$ and quilshon for every $\delta=\cf(\delta)>\aleph_0$. It follows from Theorem \ref{mt2} that there are no Magidor cardinals in this model, despite the fact that there are (at least) $\kappa$-many measurable cardinals. Moreover, ${\rm V}_\kappa$ is a model of ZFC in which there is a class of measurable cardinals and no Magidor cardinal.

\hfill \qedref{ccumandshimmer}

The main theorem of this section says that an $\omega$-limit of supercompat cardinals need not be a Magidor cardinal. Presumably, this idea can be exploited to produce a universe with a class of supercompact cardinals with no Magidor cardinal, and also to reach beyond supercompactness, see the remarks below.

Ahead of the proof we shall phrase a simple lemma, which seems useful also for other assertions of the same type. Given a function $f\in\prod\limits_{i<\delta}\theta_i$, the support of $f$ is ${\rm supp}(f)=\{i<\delta:f(i)\neq 0\}$. We say that $f$
has Easton support iff $|{\rm supp}(f)\cap\sigma|<\sigma$ whenever $\sigma$ is an inaccessible cardinal.

\begin{lemma}
\label{llll} Assume that:
\begin{enumerate}
\item[$(a)$] $\jmath:{\rm V}\rightarrow M$ is an elementary embedding and $\kappa=\crit(\jmath)$.
\item[$(b)$] The cofinality of all the cardinals in the product $\prod\limits_{i<\delta}\theta_i$ is at least $\kappa$.
\item[$(c)$] $\mu$ is an $M$-regular cardinal above $\jmath(\kappa)$.
\item[$(d)$] $\beta={\rm sup}\{\jmath f(\mu):f$ has Easton support in $\prod\limits_{i<\delta}\theta_i\}$.
\end{enumerate}
Then $\cf^M(\beta)\geq\kappa$.
\end{lemma}

\par\noindent\emph{Proof}. \newline
Assume towards contradiction that $\theta=\cf(\beta)<\kappa$. Choose $\langle f_j:j<\theta\rangle$, each $f_j$ is an Easton support function in $\prod\limits_{i<\delta}\theta_i$, such that $\beta={\rm sup}\{\jmath f_j(\mu):j<\theta\}$. Define $g={\rm sup}\{f_j:j<\theta\}+1$ and deduce from $(b)$ that $g\in\prod\limits_{i<\delta}\theta_i$ and has an Easton support. However, $\jmath g(\mu)>{\rm sup}\{f_j(\mu):j<\theta\}=\beta$, a contradiction to the very definition of $\beta$, so we are done.

\hfill \qedref{llll}

\begin{theorem}
\label{ssupercccomp} Quilshon and supercompact cardinals. \newline
It is consistent that an $\omega$-limit of supercompact cardinals is not a Magidor cardinal.
\end{theorem}

\par\noindent\emph{Proof}. \newline
Let $\langle\kappa_n:n\in\omega\rangle$ be a strictly increasing sequence of cardinals so that $\kappa_0=\aleph_0$ and $\kappa_n$ is a supercompact cardinal for every $0<n\in\omega$. We may assume that $\kappa_n$ is Laver-indestructible for every $0<n\in\omega$. Let $\lambda$ be $\bigcup\limits_{n\in\omega}\kappa_n$ and denote the set $S^\lambda_\omega$ by $S$.

We wish to define a forcing notion $\mathbb{Q}$ such that:
\begin{enumerate}
\item [$(\aleph)$] $\kappa_n$ is a supercompact cardinal after forcing with $\mathbb{Q}$, for every $0<n\in\omega$.
\item [$(\beth)$] For every $n\in\omega$ there is a regular cardinal $\kappa_n<\theta_n<\kappa_{n+1}$ for which $\pitchfork_{S,\theta_n}$ holds after forcing with $\mathbb{Q}$.
\end{enumerate}

From $(\aleph)$ we infer that $\lambda$ is still a limit of supercompact cardinals in the generic extension by $\mathbb{Q}$. From $(\beth)$ we conclude, by the quilshon claim \ref{qqqilshonfromsaure}, that $\lambda$ is not a Magidor cardinal after forcing with $\mathbb{Q}$. Let us describe this forcing notion. We indicate that $\mathbb{Q}$ may collapse some cardinals, but each Mahlo cardinal along the iteration will be preserved.

We shall define $\mathbb{Q}$ as a product of the form $\prod\limits_{n\in\omega}\mathbb{R}_n$. Each component $\mathbb{R}_n$ would be a forcing notion which adds a partial global square for the cofinality $\theta_n$, modelled basically after the forcing of Jensen. We shall see that $\mathbb{R}_n$ is $\theta_n$-closed, and we shall prove that it possesses enough completeness and enough chain condition properties in order to preserve enough cardinals. By Claim \ref{qqqilshonfromsaure} we shall get $\pitchfork_{S,\theta_n}$ at every $n\in\omega$ after forcing with $\mathbb{R}_n$. Finally, we shall see that each $\mathbb{R}_n$ preserves the supercompactness of the $\kappa_n$-s.

Assume we have accomplished the construction of $\mathbb{R}_n$ for every $n\in\omega$. Let $\mathbb{Q}=\prod\limits_{n\in\omega}\mathbb{R}_n$ with full support. Fix any $n\in\omega$ and decompose the product into the left part $\prod\limits_{m\leq n}\mathbb{R}_n$ and the right part $\prod\limits_{m> n}\mathbb{R}_n$. By the $\theta_{n+1}$-completeness of each $\mathbb{R}_m$ when $m>n$ we know that $\prod\limits_{m> n}\mathbb{R}_n$ is $\kappa_n$-complete and even $\kappa_n$-directed-closed. By the properties of $\mathbb{R}_m$ (to be proved below) we know that the supercompactness of each $\kappa_n$ is preserved by $\mathbb{R}_m$ for every $m\leq n$ and hence also by $\prod\limits_{m\leq n}\mathbb{R}_n$. It follows that $\mathbb{Q}$ preserves the supercompactness of every $\kappa_n$, so $\lambda$ is a limit of supercompact cardinals in the generic extension by $\mathbb{Q}$. However, $\lambda$ is not a Magidor cardinal after forcing with $\mathbb{Q}$, by virtue of $\pitchfork_{S,\theta_n}$ at every $n\in\omega$, so we are done.

For accomplishing the proof we have to define $\mathbb{R}_n$ for every $n\in\omega$, and to prove the asserted properties of this forcing notion.
Fix a natural number $n$ and the associated regular cardinal $\theta_n\in(\kappa_n,\kappa_{n+1})$.
Let $\mathbb{R}=\mathbb{R}_n$ be an iteration with reverse Easton support over all the regular cardinals to add a partial global square at the cofinality $\theta_n$. It means that we take direct limits at inaccessible stages and inverse limits at other limit stages. In the successor stage, let $\mathbb{R}_{\beta+1}=\mathbb{R}_\beta\ast\mathbb{S}_\beta$ when $\mathbb{S}_\beta=\{\emptyset\}$ if $\beta$ is not a regular cardinal, and we define $\mathbb{S}_\beta$ for a regular $\beta$ as follows.

A condition $p\in\mathbb{S}_\beta$ is an approximation to a global square sequence $p=\langle C_\gamma:\gamma\leq\delta<\beta, {\rm otp}(C_\gamma)\leq\theta_n\rangle$ such that:
\begin{enumerate}
\item [$(a)$] $C_\gamma=\emptyset$ or $C_\gamma$ is a club in $\gamma$.
\item [$(b)$] If $\gamma>\cf(\gamma)=\theta_n$ then $C_\gamma$ is a club in $\gamma$.
\item [$(c)$] The sequence is coherent, i.e. if $\xi\in{\rm acc}(C_\gamma)$ then $C_\gamma\cap\xi=C_\xi$.
\end{enumerate}
For the order, assume $p,q\in\mathbb{S}_\beta$ and let $p\leq_{\mathbb{S}_\beta}q$ iff $\delta_p\leq\delta_q$ and $C^p_\gamma=C^q_\gamma$ for every $\gamma\leq\delta_p$.

The component $\mathbb{S}_\beta$ of the forcing $\mathbb{R}$ is $\theta_n$-complete. Indeed, assume $\zeta<\theta_n$ and $\langle p_i:i<\zeta\rangle$ is an increasing sequence of conditions in $\mathbb{S}_\beta$. Let $q$ be the union $\bigcup\limits_{i<\zeta}p_i$ with $C_{\zeta+1}=\emptyset$ appended as a top element. Notice that $\cf(\zeta)\leq\zeta<\theta_n$ so we may add the empty set as a last element. It follows that $q\in\mathbb{S}_\beta$ and $p_i\leq q$ for every $i<\zeta$.

Recall that if $\mathbb{P}$ is a forcing notion, $p\in\mathbb{P}$ and $\mu$ is an infinite cardinal then $\Game_\mu(p,\mathbb{P})$ is a two-players game which lasts $\mu$ moves. In the $\alpha$-th move, player I tries to choose $p_\alpha\in\mathbb{P}$ such that $p\leq p_\alpha$ and $\beta<\alpha\Rightarrow q_\beta\leq p_\alpha$. Player II tries to choose $q_\alpha\in\mathbb{P}$ such that $p_\alpha\leq q_\alpha$. Player I wins iff she has a legal move for every $\alpha<\mu$. The forcing $\mathbb{P}$ is $<\mu$-strategically complete iff player I has a winning strategy in $\Game_\mu(p,\mathbb{P})$ for every $p\in\mathbb{P}$.

We claim now that $\mathbb{S}_\beta$ is $<\beta$-strategically closed. For proving it, assume $\mu<\beta$ and we have the usual two-players game of length $\mu$. Let $\mathcal{D}=\langle D_\delta:\delta<\mu\rangle$ be a partial square sequence along the ordinals of $S^\mu_{\theta_n}$. It exists, due to the induction hypothesis. The strategy of the good player will be to choose at every limit stage $\eta$ the set $\{\xi_i:i\in D_\eta\}$, stipulating $\xi_i=\ell g(p_i)$ (meaning the length of $p_i$). It follows that this is a winning strategy.

It remains to show that $\mathbb{R}=\mathbb{R}_n$ preserves the supercompactness of $\kappa_n$. Denote $\kappa_n$ by $\kappa$, and fix any cardinal $\mu$. We shall prove that if $\kappa$ is $2^{(\mu^{<\kappa})}$-supercompact then it remains $\mu$-supercompact after forcing with $\mathbb{R}$. In particular, if $\kappa$ is supercompact then forcing with $\mathbb{R}$ preserves it.

For the fixed $\mu$ let $\tau=2^{(\mu^{<\kappa})}$, and let $\jmath:{\rm V}\rightarrow M$ be any $\tau$-supercompact embedding for which $\kappa=\crit(\jmath)$ and $\jmath(\kappa)>\tau$. By Silver's criterion we are looking for an $M$-generic $H$ with respect to the forcing notion $\jmath(\mathbb{R})$ such that $\jmath``G\subseteq H$. If we can create such an $H$ then we will be able to extend $\jmath$ into $\jmath^+:{\rm V}[G]\rightarrow M[H]$. The lower part of $H$ is inherited from $G$ by the Easton support (as described below), and the rest will be constructed by defining the pertinent strong master condition.

Notice that forcing with $\mathbb{R}$ after the $\tau$-stage adds no sets of size $\tau$ or less. Consequently, this part of $\mathbb{R}$ neither adds nor destroys supercompact measures over $\mathcal{P}_\kappa\mu$. Hence we may concentrate on the first $\tau$ stages of $\mathbb{R}$, which we still call $\mathbb{R}$ for simplicity.

As a first step let $H\upharpoonright\kappa=G\upharpoonright\kappa$. This is justified by the Easton support. Indeed, for every condition $p\in\mathbb{R}$ we have ${\rm supp}(p)\cap\kappa$ is bounded, hence ${\rm supp}(\jmath(p)\cap\jmath(\kappa)) = \jmath({\rm supp}(p)\cap\kappa) = {\rm supp}(p)\cap\kappa$. Moreover, $M$ agrees with ${\rm V}$ up to the stage $\tau$ on $\mathbb{R}$, so we can define $H\upharpoonright\tau=G$.

Define $m=\bigcup\{\jmath(p)\upharpoonright[\jmath(\kappa),\jmath(\tau)]:p\in G\}$, meaning that we take at every point the concatenation of all the square sequences at this coordinate.
We have to show that $m$ is a name of a condition in $\jmath(\mathbb{R})\upharpoonright[\jmath(\kappa),\jmath(\tau)]$. Since ${}^\tau M\subseteq M$ (by the choice of $\jmath$) we know that $\jmath``\mathbb{R}\in M$. It follows that $m\in M[H\upharpoonright\tau]$. Likewise, if $\cf^M(\mu)=\mu$ then $\cf^M(\ell g (\bigcup\limits_{p\in G}\jmath(p)(\mu)))\geq \kappa$ due to Lemma \ref{llll} so we can set $C_{\mu+1}=\emptyset$ without violating item $(b)$ in the definition of the members of $\mathbb{R}$, and thus complete $m$ to a condition.

Finally, $m$ is a strong master condition by the construction. More precisely, $m$ has been defined on the interval $[\jmath(\kappa),\jmath(\tau)]$ in such a way that exceeds every embedding of conditions in $G$, and we may have to fix $m$ below $\jmath(\kappa)$ in order to get a strong master condition.
We can force now over ${\rm V}[G]$ the existence of a ${\rm V}$-generic set $H$ for $\jmath(\mathbb{R})$ such that $m\in H$. Since ${}^\tau M\subseteq M$ we know that $H$ is also $M$-generic.

We have defined $H$ and proved that $\jmath``G\subseteq H$. Denote $K=H\upharpoonright[\tau^+,\jmath(\tau)]$ and notice that $\jmath^+:{\rm V}[G]\rightarrow M[H]$ is definable in ${\rm V}[G][K]$. In particular, there exists in ${\rm V}[G][K]$ a normal ultrafilter over $\mathcal{P}_\kappa\mu$. However, the part of $\mathbb{R}$ which is added by $K$ over ${\rm V}[G]$ does not add sets of size less than $\tau^+$, and hence the normal measure on $\mathcal{P}_\kappa\mu$ belongs to ${\rm V}[G]$, so we are done.

\hfill \qedref{ssupercccomp}

\begin{remark}
\label{r5}
Unlike the forcing in Corollary \ref{ccumandshimmer}, which adds a quilshon for every $\delta=\cf(\delta)>\aleph_0$, in the above theorem we add a quilshon only for an unbounded set of regular cardinals below $\lambda$. Essentially, the reason is the strong reflection properties of the supercompact cardinal.

The result is that we can make sure that $\lambda$ is not Magidor, but maybe there are Magidor cardinals in the universe (even below $\lambda$). However, we believe that it is consistent to have a class of supercompact cardinals with no Magidor cardinals at all.
\end{remark}

\hfill \qedref{r5}

The process of adding an unbounded sequence of partial squares below $\lambda$ in order to eliminate its Magidority can be applied also in the context of rank-to-rank embeddings. We shall prove that an $\omega$-limit of Magidor cardinals can be a non-Magidor cardinal. Of course, such a theorem assumes that there are at least $\omega$ Magidor cardinals. Moreover, we will have to assume that they come from instances of I1, and that the critical points satisfy some preliminary requirement. The main reason is that we must preserve the Magidority of the members in the $\omega$-sequence of Magidor cardinals. It follows that we can preserve I1 (and consequently, the pertinent Magidority). The proof is based on the ideas of the former proof (for supercompact cardinals), and in some sense it is a bit simpler.

\begin{theorem}
\label{llimitofmagidor} Non-limitude of Magidority. \newline
Assume that $I1(\kappa_n,\lambda_n)$ holds, and $\lambda_n<\kappa_{n+1}$ for every $n\in\omega$. \newline
Let $\tau$ be $\bigcup\limits_{n\in\omega}\lambda_n$. \newline
There is a forcing notion $\mathbb{Q}$ such that if $G_{\mathbb{Q}}\subseteq\mathbb{Q}$ is generic then in ${\rm V}[G_{\mathbb{Q}}]$ every $\lambda_n$ is a Magidor cardinal but $\tau$ is not a Magidor cardinal.
\end{theorem}

\par\noindent\emph{Proof}. \newline
We choose a regular uncountable cardinal $\theta_0<\kappa_0$, and for every $0<n\in\omega$ we choose a regular cardinal $\theta_n$ such that $\lambda_{n-1}<\theta_n<\kappa_n$. We define the forcing notion $\mathbb{R}_n$ to be the forcing which adds a partial square for the cofinality $\theta_n$. Define $\mathbb{Q}=\prod\limits_{n\in\omega}\mathbb{R}_n$ with full support, and let $G_{\mathbb{Q}}\subseteq\mathbb{Q}$ be a generic set.

Fix a natural number $n\in\omega$, and split $\mathbb{Q}$ into a product of $\prod\limits_{m\leq n}\mathbb{R}_n$ and$\prod\limits_{m> n}\mathbb{R}_n$. The upper part would be $\theta_m$-complete, so $I1(\kappa_n,\lambda_n)$ is preserved. Indeed, no new subsets of ${\rm V}_{\lambda_n+1}$ are introduced by this part of the forcing. We shall see below that $\mathbb{R}_\ell$ preserves $I1(\kappa_n,\lambda_n)$ for every $\ell\leq n$, and hence the lower part $\prod\limits_{m\leq n}\mathbb{R}_n$ preserves $I1(\kappa_n,\lambda_n)$ as well. This is the main burden of the proof, and if we succeed then we conclude with $I1(\kappa_n,\lambda_n)$ for every $n\in\omega$ in ${\rm V}[G_{\mathbb{Q}}]$. On the other hand, forcing with $\mathbb{Q}$ adds a partial square for every $\theta_n$. Since these cardinals are unbounded in $\lambda$ we infer that $\lambda$ is not a Magidor cardinal in ${\rm V}[G_{\mathbb{Q}}]$.

As usual, the part of $\mathbb{R}_n$ above $\lambda_n^+$ cannot add new sets of cardinality less than $\lambda_n^+$ (by the distributivity of this component), and hence has no influence on $I1(\kappa_n,\lambda_n)$. We focus on $\mathbb{R}_n\upharpoonright\lambda_n^+$ which will be denoted by $\mathbb{R}$. We denote the generic set that we choose for $\mathbb{R}$ by $G$, and let $\kappa_n=\kappa, \lambda_n=\lambda$.

The key point in proving that forcing with $\mathbb{R}$ preserves $I1(\kappa_n,\lambda_n)$ is that the conditions in $\mathbb{R}$ are elements of ${\rm V}_{\lambda+1}$. Indeed, if $p\in\mathbb{R}$ then ${\rm dom}(p)\cap\jmath^n(\kappa)$ is bounded for every $n\in\omega$ as we use Easton support. In particular, $p\in{\rm V}_{\lambda+1}$ and hence $p\in{\rm dom}(\jmath)$. Moreover, every set in ${\rm V}_{\lambda+1}[G]$ has a name in ${\rm V}_{\lambda+1}[G]$ (more precisely, a $\lambda$-sequence of members of ${\rm V}_{\lambda+1}\times {\rm V}_\lambda$ which can be coded as an element of ${\rm V}_{\lambda+1}[G]$).

Aiming to lift $\jmath:{\rm V}_{\lambda+1}\rightarrow{\rm V}_{\lambda+1}$ into $\jmath^+:{\rm V}_{\lambda+1}[G]\rightarrow {\rm V}_{\lambda+1}[H]$ we must define a suitable generic set $H$.
We emphasize that in this case our purpose is to get ${\rm V}_{\lambda+1}[G]={\rm V}_{\lambda+1}[H]$, since we are looking for an embedding from ${\rm V}_{\lambda+1}[G]$ to itself in the generic extension.
Observe that if $p\in\mathbb{R}$ then $\jmath(p)\in\mathbb{R}$ by the definition of $\mathbb{R}$ and the elementarity of $\jmath$, though $p\neq\jmath(p)$ in general. Hence the set $H$ should be $\mathbb{R}$-generic over ${\rm V}_{\lambda+1}$. We also must make sure that $\jmath``G\subseteq H$ so we define a master condition and depict a generic set $H$ which contains it.

We commence with a description of a master condition $m\in\mathbb{R}$. The condition $m$ will not be a strong master condition (i.e., $m\geq\jmath(p)$ for every $p\in G$) but it will satisfy $p\upharpoonright\kappa\cup m\geq\jmath(p)$ for every $p\in G$ which is enough in order to use Silver's criterion. The advantage is that we can concentrate on conditions $p$ for which ${\rm dom}(p)\cap\kappa=\emptyset$. Set:
$$
m=\bigcup\{\jmath(p):p\in {\rm V}_{\lambda+1}\cap \name{G}, {\rm dom}(p)\cap\kappa=\emptyset\}.
$$
By the union we mean that we take concatenations of the square approximation sequences at each stage, and then we take all the sequences as our condition. Of course, we have to show that this definition of $m$ is actually a condition in $\mathbb{R}$. For this end, we shall prove that $m$ has Easton support and that the pertinent cofinality allows us to add the empty set as a last element in the $\mu$-th coordinate for every $\mu=\cf^M(\mu)$ above $\kappa$.

Assume, firstly, that $\mu$ is a regular cardinal in $M$ and $\kappa\leq\mu<\lambda$, so $\jmath^n(\kappa)\leq\mu<\jmath^{n+1}(\kappa)$ for some $n\in\omega$. Set $\name{\beta}={\rm sup}\{{\rm dom}\jmath(p)(\mu):p\in\name{G}\}$, assuming without loss of generality that this fact is forced by the empty condition so we can treat $\name{\beta}$ as an ordinal $\beta$. Let $A$ be the set of all functions with Easton support in $\prod\{\eta:\jmath^{n-1}(\kappa)\leq\eta<\jmath^n(\kappa), \eta=\cf(\eta)\}$. By the reasoning of Lemma \ref{llll}, using genericity arguments, we infer that $\cf^M(\beta)\geq\kappa$. Since $\theta_n<\kappa_n=\kappa\leq\cf^M(\beta)$ we can set $C_{\mu+1}=\emptyset$, so the definition of $m$ as a condition can be accomplished.

Let us try to show that $m$ has Easton support. Recall that ${\rm dom}(p)\cap[\jmath^n(\kappa),\jmath^{n+1}(\kappa))$ is bounded in $\jmath^{n+1}(\kappa)$ for every $n\in\omega$ and every $p\in\mathbb{R}$. Consequently, if $\jmath^n(\kappa)<\mu=\cf(\mu)\leq\jmath^{n+1}(\kappa)$ then ${\rm dom}(m)\cap\mu$ is the union of at most $\jmath^n(\kappa)$ many sets (since $|\mathbb{R}\upharpoonright \jmath^n(\kappa)|=\jmath^n(\kappa)$), each of which is bounded in $\mu$. It follows that $|{\rm dom}(m)\cap\mu|\leq\jmath^n(\kappa)<\mu$ as desired.

Finally, let $p\in G$ be a condition for which ${\rm dom}(p)\cap\kappa=\emptyset$. We need showing that $\jmath(p)\leq m$, and we shall prove that $p\upharpoonright\alpha\Vdash\jmath(p)(\alpha)\leq m(\alpha)$ for every $\alpha$. Recall that if $m\geq\jmath(p)$ whenever ${\rm dom}(p)\cap\kappa=\emptyset$ then $p\upharpoonright\kappa\cup m\geq\jmath(p)$ for every $p\in G$, which is sufficient.

There is nothing to worry about for $\alpha<\kappa$ by the assumption that ${\rm dom}(p)\cap\kappa=\emptyset$. Likewise, there is nothing to worry about for $\kappa\leq\alpha<\jmath(\kappa)$ since $\jmath(p)(\alpha)$ is empty in these cases as $\kappa=\crit(\jmath)$ (so the first element of ${\rm dom}(p)$ will be sent by $\jmath$ at least to $\jmath(\kappa)$). Suppose $\alpha\geq\jmath(\kappa)$. By definition, $m(\alpha)$ is a name for an upper bound of $\{\jmath(q)(\alpha):q\in \name{G}\}$. In particular, $p\upharpoonright\alpha\Vdash m(\alpha)\geq\jmath(p\upharpoonright\alpha)(\alpha)=\jmath(p)(\alpha)$.

Having established the definition of $m$ and the fact that it can serve as a master condition, we can choose a generic set $G$ such that $m\in G$. We define $H=G$, so $H$ is $\mathbb{R}$-generic and includes $m$. By Silver's criterion we can lift our embedding $\jmath$ into the extension $\jmath^+:{\rm V}_{\lambda+1}[G]\rightarrow {\rm V}_{\lambda+1}[H] = {\rm V}_{\lambda+1}[G]$ and conclude that $I1(\kappa,\lambda)=I1(\kappa_n,\lambda_n)$ is preserved by $\mathbb{R}$ as required.

\hfill \qedref{llimitofmagidor}

\begin{remark}
\label{belowi1} The non-Magidority of an $\omega$-limit of Magidor cardinals is not confined to a limit of I1 instances. Starting from the above setting, we can add $\lambda^+$ Cohen subsets of $\aleph_1$ after the forcing of the above theorem. The result is a non-Magidor limit of Magidor cardinals below $2^{\aleph_1}$, and in particular they are not I1.
\end{remark}

\hfill \qedref{belowi1}

The above theorem yields similar consistency results for other notions of large cardinals, and we phrase below a typical one.
Recall that a Magidor cardinal can be very small in the $\beth$-scale, by Claim \ref{ssmallmagidor}. However, the above theorem shows that the natural way to construct a Magidor cardinal as a limit of large cardinals, fails even under very strong hypotheses.
This kind of statements creates the impression that a Magidor cardinal is a very large cardinal, at least from this point of view.
The moral of the corollary below is that sometimes it doesn't help to stand on the shoulders of giants:

\begin{corollary}
\label{hhugelimit} A huge corollary. \newline
It is consistent (under the assumptions of the previous theorem) that an $\omega$-limit of huge cardinals is not a Magidor cardinal.
The same holds if we replace huge by $n$-huge for every $n\in\omega$.
\end{corollary}

\par\noindent\emph{Proof}. \newline
If $I1(\kappa,\lambda)$ then there is a normal ultrafilter $U$ over $\kappa$ which concentrates on the set $\{\delta<\kappa:\delta$ is $n$-huge for every $n\in\omega\}$. Actually, such $U$ exists even under the weaker assumption $I3(\kappa,\lambda)$. By Theorem \ref{llimitofmagidor}, choose such $\delta_m$ for every $m\in\omega$, which satisfies $\lambda_m<\delta_m<\kappa_{m+1}$ in ${\rm V}[G]$. It follows that $\lambda=\bigcup\limits_{m\in\omega}\lambda_m = \bigcup\limits_{m\in\omega}\delta_m$ is a non-Magidor cardinal, limit of cardinals which are $n$-huge for every $n\in\omega$.

\hfill \qedref{hhugelimit}

The assumption that the intervals $[\kappa_n,\lambda_n)$ are pairwise disjoint (both in Theorem \ref{llimitofmagidor} and in Corollary \ref{hhugelimit}) is quite reasonable. In fact, if there are infinitely many instances of $I1(\kappa_n,\lambda_n)$ for which $m\neq n \Rightarrow \kappa_m\neq \kappa_n$ then there is also a sequence of pairs which satisfy the disjointness assumption (although not necessarily a sub-sequence of the original one).

We comment that an interlaced sequence of the form $\kappa_n<\kappa_{n+1}<\lambda_n<\lambda_{n+1}$ and $I1(\kappa_n,\lambda_n)$ for every $n\in\omega$ might prove as the correct assumption for showing the Magidority of $\lambda=\bigcup\limits_{n\in\omega}\lambda_n$.
More generally, if $\kappa$ is $2$-super huge then partial squares and even quilshon fail for every $\delta=\cf(\delta)\geq\kappa$. We also mention Martin's maximum, from which $\neg(\pitchfork_{\lambda,\aleph_2})$ follows. If this can be generalized to higher cofinalities, then we might be able to prove positive Magidority results.

\newpage

\section{Open problems}

We list a few open problems concerning Magidor cardinals and the related cardinal $\alpha_M$. For brevity, $\lambda$ is a Magidor cardinal throughout this section. Some of the problems have an analog for J\'onsson cardinals, including the first one.

\begin{question}
\label{q1} Is it consistent that $\aleph_\omega$ is a Magidor cardinal?
\end{question}

A related problem can be asked about the limit cardinals below any Magidor cardinal. We have shown the consistency of $\lambda<2^{\aleph_1}$, so a Magidor cardinal need not be a limit of strongly inaccessible cardinals. However, we may ask:

\begin{question}
\label{q2} Is it consistent that a Magidor cardinal is not a limit of weakly inaccessible cardinals?
\end{question}

A negative answer settles, of course, the first question. Another way to attack the first question sends us to the realm of determinacy.
Recall that AD implies that $\aleph_\omega$ is a Rowbottom cardinal, by \cite{MR0469769}. By and large, AD has a salient propensity to combinatorial assertions with infinite exponent.
The most conspicuous assertion is $\omega_1\rightarrow(\omega_1)^{\omega_1}$, and it follows from this property that $\aleph_1$ is a Magidor cardinal. The case of $\aleph_\omega$ invites a similar conclusion, so we ask:

\begin{question}
\label{q25} Assume AD. Is it provable that $\aleph_\omega$ is a Magidor cardinal?
\end{question}

We recall that a limit of J\'onsson cardinals need not be a J\'onsson cardinal, by a result of Kunen. Similarly, under the assumption that there are many instances of I1 we have shown that the limit can be a non-Magidor cardinal (see Theorem \ref{llimitofmagidor}, under an additional assumption on the sequence). The I1 assumption was essential in the proof, but non limitude of Magidority is consistent even without having I1, see remark \ref{belowi1}.

The other side of the coin is the axiom I0, which says that there exists an elementary embedding $\jmath:L({\rm V}_{\lambda+1}) \rightarrow L({\rm V}_{\lambda+1})$ such that $\crit(\jmath)<\lambda$. It has been proved by Laver that I0 is strictly stronger than I1. Presumably, the forcing in Theorem \ref{llimitofmagidor} destroys the property of being I0.

\begin{question}
\label{q30} Let $\lambda$ be an $\omega$-limit of Magidor cardinals which are I0. Is it consistent that $\lambda$ is not a Magidor cardinal?
\end{question}

It is known that if $I1(\kappa,\lambda)$ then there are $\kappa_0<\kappa,\lambda_0<\lambda$ for which $I2(\kappa_0,\lambda_0)$. Consequently, the first Magidor cardinal is not I1. The following is natural:

\begin{question}
\label{q31} Is it consistent that the first Magidor cardinal is I2?
\end{question}

Finally, it is well known that the first J\'onsson cardinal is either weakly inaccessible or a singular cardinal with countable cofinality.
It follows from Corollary \ref{bbbb} that the first J\'onsson cardinal need not be the first Magidor cardinal.
We may ask if the first Magidor cardinal has some reflection properties with respect to J\'onsson cardinals:

\begin{question}
\label{q33} Is it consistent that the first J\'onsson cardinal is also the first Magidor cardinal?
\end{question}

Let us turn to the relation between Magidor cardinals and other large cardinals. In the former sections we have seen on the one hand that Magidority is a consequence of I1 and I2, and on the other hand that a global non-Magidority is consistent with the existence of a class of measurable cardinals. We may ask about larger cardinals, and the following two-fold question is typical:

\begin{question}
\label{q4} Is it consistent that there is one supercompact cardinal and no Magidor cardinal? Is it consistent that there is a class of supercompact cardinals and no Magidor cardinal?
\end{question}

In fact, one should state this question with respect to smaller notions (e.g., strong compactness). It begins at the point of the failure of square (but the amount of square is essential here). Anyhow, Theorem \ref{ssupercccomp} is suggestive, and we believe that the answer is positive. From the other direction of large cardinals:

\begin{question}
\label{q5} Assume $\lambda$ is I3 (i.e. there exists an elementary embedding $\jmath:{\rm V}_\lambda\rightarrow {\rm V}_\lambda$). Does it follow that $\lambda$ is a Magidor cardinal?
\end{question}

We comment that I3 seems to be fairly different from I1 and I2, so one may try the opposite direction (i.e., proving the consistency of I3 with non-existence of Magidor cardinals). One way to do it is to force weak versions of square or quilshon while preserving stronger assumptions of large cardinals. We mention here the huge cardinals. The partial square used in Theorem \ref{ssupercccomp} becomes a bit more problematic in the huge environment, and much more problematic above a super huge cardinal.

By Corollary \ref{hhugelimit}, an $\omega$-limit of huge cardinals can be non-Magidor.
Moreover, hugeness can be sharpened to $n$-hugeness for every $n\in\omega$. However, it is important to bound the target of the elementary embedding.
Therefore, we may ask:

\begin{question}
\label{q6} Is it consistent that a non-Magidor cardinal $\lambda$ is a limit of super huge cardinals?
Does there exist a strong enough notion of large cardinals which ensures that an $\omega$-limit of its type is a Magidor cardinal?
\end{question}

From the other direction, we have seen that each Magidor cardinal is J\'onsson. Likewise, if $\lambda$ is Magidor then $\lambda>2^{\aleph_0}$ and $\omega$-closed. Recall that every Rowbottom cardinal above $2^{\aleph_0}$ is $\omega$-closed (see \cite{MR654576}). A Magidor cardinal need not be Rowbottom (as commented by the referee of the paper, by forcing $\lambda$ to be below $2^{\aleph_1}$), but we raise the following:

\begin{question}
\label{qrowbottom} Is it true that every Magidor cardinal is $\alpha_M$-Rowbottom?
\end{question}

The following problems are connected with $\alpha_M$, the first ordinal for which the coloring omits colors.
We have proved that if $\theta<\alpha_M$ and $\cf(\theta)>\omega$ then $\theta^\omega<\alpha_M$. The motivation is that if we proved that $\alpha_M$ is $\omega$-closed then we could prove the consistency of $\alpha_M=\aleph_2$ (which is the lower bound on $\alpha_M$). However, the proviso $\cf(\theta)>\omega$ seems obstinate:

\begin{question}
\label{q7} Is it consistent that $\lambda$ is a Magidor cardinal and $\alpha_M=\theta^+$ when $\cf(\theta)=\omega$? Is it possible for a Magidor cardinal which comes from $I1(\kappa,\lambda)$?
\end{question}

Concerning this problem, one may wish to focus on the first Magidor cardinal, in which case $\theta<\alpha_M<\lambda$ and hence $\theta$ itself is not a Magidor cardinal. This assumption gives some hope to imitate the proof of the uncountable cofinality, but there are still some obstacles. Anyway, we may ask more generally:

\begin{question}
\label{q8} Is it consistent that $\lambda$ is a Magidor cardinal and $\alpha_M$ is a successor of a singular cardinal?
\end{question}

We indicate that pcf arguments may be helpful to refute the above possibility (by introducing $\alpha_M$ as the true cofinality of some sequence of cardinals for which there are no Magidor functions). Another possible direction is to begin with $\alpha_M$ as a successor of a large cardinal and to try to singularize this cardinal, so:

\begin{question}
\label{q9} Is it consistent that $\lambda$ is a Magidor cardinal and $\alpha_M$ is a successor of a measurable cardinal? Is it consistent that $\alpha_M$ is a measurable cardinal?
\end{question}

A general question is the preservation of Magidority under forcing extensions.
In particular, we should investigate the influence of Prikry-type forcing notions on Magidority:

\begin{question}
\label{q11} Suppose $\lambda$ is a Magidor cardinal, limit of measurable cardinals.
Let $\mathbb{P}$ be the diagonal Prikry forcing with respect to $\lambda$. Is it true that $\lambda$ is still a Magidor cardinal after forcing with $\mathbb{P}$?
\end{question}

Notice that the assumption on $\lambda$ in the above question holds when $\lambda$ comes from an instance of I1.

It has been proved for $\alpha>\cf(\alpha)>\aleph_0$ that $\lambda\rightarrow [\lambda]^{\aleph_0\text{-bd}}_\alpha$ implies $\lambda\rightarrow [\lambda]^{\aleph_0\text{-bd}}_{\alpha,<\alpha}$. We encountered the recurrent problem of the countable cofinality case, which seems to be essential for several reasons when dealing with Magidor cardinals. We may ask:

\begin{question}
\label{q12} Assume $\lambda\rightarrow [\lambda]^{\aleph_0\text{-bd}}_\alpha$, and $\cf(\alpha)=\aleph_0$. Is it provable that $\lambda\rightarrow [\lambda]^{\aleph_0\text{-bd}}_{\alpha,<\alpha}$?
\end{question}

The parallel for J\'onsson cardinals holds true, by assuming that the number of colors is not a J\'onsson cardinal, due to Kleinberg. It means that $\lambda\rightarrow [\lambda]^{<\omega}_{\alpha,<\alpha}$ whenever $\alpha$ is not a J\'onsson cardinal. We mention, en route, that $\lambda\rightarrow [\lambda]^{<\omega}_{\alpha,<\alpha}$ is consistent even if $\alpha$ is J\'onsson, assuming e.g. that $\lambda$ carries a Rowbottom filter. But this assumption, as it is, cannot be generalized to Magidority.

Unlike J\'onssonicity, for which the crux of the matter is whether $\alpha$ is J\'onsson or not, the limitation in our context is whether $\alpha$ has countable cofinality. We may wish to separate between a Magidor $\alpha$ and a non-Magidor $\alpha$ in the above problem. The assumption that $\alpha$ is not a Magidor cardinal holds always for the first Magidor cardinal $\lambda$.

\begin{question}
\label{q125} Let $\lambda_0<\lambda_1$ be Magidor cardinals, and $\alpha_M^0,\alpha_M^1$ the associated cardinals. Is it consistent that $\alpha_M^1<\alpha_M^0$? What about a Magidor cardinal, limit of the sequence of Magidor cardinals $\langle\lambda_n:n\in\omega\rangle$ assuming that $\langle\alpha_M^n:n\in\omega\rangle$ is strictly increasing?
\end{question}

Concerning the second part, it might be helpful to distinguish two cases. Let $\alpha=\bigcup\limits_{n\in\omega}\alpha_M^n$. If $\alpha=\lambda$ and it is a Magidor cardinal then necessarily $\alpha_M<\alpha_M^n$ for almost every $n\in\omega$.

We conclude with a problem about the relation between $\alpha_J$ and $\alpha_M$. We have seen that $\alpha_J\leq\alpha_M$ for every Magidor cardinal, and if $I1(\kappa,\lambda)$ then it seems that one can force strict inequality by adding $\alpha_J$-many Cohen subsets to $\aleph_0$ (recall that $2^{\aleph_0}<\alpha_M$). It seems harder to prove equality:

\begin{question}
\label{q13} Is it consistent that $\alpha_J=\alpha_M$? Is it consistent for a Magidor cardinal which comes from I1?
\end{question}

We include here a partial answer to Question \ref{q12}, proved by the referee of the paper:

\begin{proposition}
\label{solreferee} Assume $\lambda\rightarrow [\lambda]^{\aleph_0\text{-bd}}_\alpha$, and $\alpha$ is not a Magidor cardinal. Then $\lambda\rightarrow [\lambda]^{\aleph_0\text{-bd}}_{\alpha,<\alpha}$.
\end{proposition}

\par\noindent\emph{Proof}. \newline
The theorem has been proved in the case of $\cf(\alpha)>\omega$, so assume $\cf(\alpha)=\omega$ and fix an increasing sequence of regular uncountable cardinals $\langle\alpha_n:n\in\omega\rangle$ with limit $\alpha$. Let $f:[\lambda]^{\aleph_0\text{-bd}}\rightarrow\alpha$ be any coloring. We have to find a set $A\in[\lambda]^\lambda$ such that $|f``[A]^{\aleph_0\text{-bd}}|<\alpha$.

By the assumption that $\alpha$ is not Magidor, fix a function $g:[\alpha]^{\aleph_0\text{-bd}}\rightarrow\alpha$ which exmeplifies the non-Magidority of $\alpha$. Likewise, for every $\delta<\lambda$ we fix a coding reals function $r_\delta:[\delta]^\omega\rightarrow{}^\omega ([\omega_1]^\omega)$.

We define a new coloring $h:[\lambda]^{\aleph_0\text{-bd}}\rightarrow\alpha$ as follows. Given $x\in[\lambda]^{\aleph_0\text{-bd}}$ we ask whether there exists a limit ordinal $\beta<\omega_1$ such that ${\rm otp}(x)=\beta+\omega$. If not, let $h(x)=0$. If the answer is yes, let $\langle x_i:i<\beta+\omega\rangle$ be an increasing enumeration of the members of $x$ and let $\delta={\rm sup}(x)$. Set $s=r_\delta(\{x_{\beta+i}:i<\omega\})$. We can express $s$ as $\langle a_j:j<\omega\rangle$ where $a_j\in[\beta]^\omega$ for every $j\in\omega$. Denote $\{x_i:i\in a_j\}$ by $y_j$. Now if $\{f(y_j):j\in\omega\}\in[\alpha]^{\aleph_0\text{-bd}}$ then $h(x)=g(\{f(y_j):j\in\omega\})$. Otherwise, let $h(x)=0$.

Choose a set $A\in[\lambda]^\lambda$ such that $h``[A]^{\aleph_0\text{-bd}}\neq\alpha$. 

We claim that $|f``[A]^{\aleph_0\text{-bd}}|<\alpha$. Assume toward contradiction that this is not the case, i.e. $|f``[A]^{\aleph_0\text{-bd}}|=\alpha$. By induction on $n\in\omega$ we choose a set $B_n\subseteq f``[A]^{\aleph_0\text{-bd}}$ such that:
\begin{enumerate}
\item [$(\alpha)$] $|B_n|\geq\alpha_n$.
\item [$(\beta)$] ${\rm sup}(B_n)<{\rm min}(B_{n+1})<\alpha$.
\item [$(\gamma)$] $\exists\mu_n<\lambda$ such that $B_n\subseteq f``[A\cap\mu_n]^{\aleph_0\text{-bd}}$.
\end{enumerate}
How do we choose these sets? For $B_0$ we choose any bounded subset of $\alpha$ of size at least $\alpha_0$ so $(\alpha)$ and $(\beta)$ are satisfied. Requirement $(\gamma)$ can be arranged by the assumption toward contradiction. The inductive step is similar. Suppose $B_i$ has been defined for $i<n$ and $n>0$. Let $\alpha'={\rm sup}(B_{n-1})<\alpha$. By the assumption toward contradiction we have $|f``[A\cap\mu_n]^{\aleph_0\text{-bd}}\setminus\alpha'|=\alpha$. Since $\cf(\alpha)=\omega$ and $\alpha_n=\cf(\alpha_n)>\omega$ we can find $\alpha``<\alpha$ such that $\alpha'<\alpha``$ and $|(f``[A\cap\mu_n]^{\aleph_0\text{-bd}}\setminus\alpha')\cap\alpha``|\geq\alpha_n$. Set $B_n= (f``[A\cap\mu_n]^{\aleph_0\text{-bd}}\setminus\alpha')\cap\alpha``$ and verify the above requirements.

We define $B^*=\bigcup\limits_{n\in\omega}B_n$. The following properties can be derived from the construction (observe, in particular, that $(b)$ follows from $(\beta)$):
\begin{enumerate}
\item [$(a)$] $|B^*|=\alpha$.
\item [$(b)$] For every $x\in[B^*]^{\aleph_0\text{-bd}}$ there is $n\in\omega$ such that $x\subseteq\bigcup\limits_{i<n}B_i$.
\item [$(c)$] For every $n\in\omega$ there is $\mu_n<\lambda$ such that $\bigcup\limits_{i<n}B_i\subseteq f``[A\cap\mu_n]^{\aleph_0\text{-bd}}$.
\end{enumerate}
Since $|B^*|=\alpha$ we know that $g``[B^*]^{\aleph_0\text{-bd}}=\alpha$. We shall use this fact in order to show that $h``[A]^{\aleph_0\text{-bd}}=\alpha$. So fix any color $\eta<\alpha$, and we shall designate a set $x\in [A]^{\aleph_0\text{-bd}}$ for which $h(x)=\eta$.

Firstly, we choose $z\in[B^*]^{\aleph_0\text{-bd}}$ such that $g(z)=\eta$. Now we recover $x$ from $z$ as follows. We choose $n\in\omega$ for which $z\subseteq \bigcup\limits_{i<n}B_i$. For this $n$ we had $\mu_n<\lambda$ which satisfies $(c)$. Let $\mu=\mu_n$, so $\bigcup\limits_{i<n}B_i\subseteq f``[A\cap\mu]^{\aleph_0\text{-bd}}$. Enumerate the members of $z$ by $\langle\zeta_j:j\in\omega\rangle$.

For every $j\in\omega$ pick a set $y_j\in [A\cap\mu]^{\aleph_0\text{-bd}}$ such that $f(y_j)=\zeta_j$. Choose $y\in [A\cap\mu]^{\aleph_0\text{-bd}}$ such that $\mu<{\rm min}(y)$ and ${\rm otp}(\bigcup\limits_{j\in\omega}y_j\cup y)=\beta$ for some limit ordinal $\beta<\omega_1$. Denote $\bigcup\limits_{j\in\omega}y_j\cup y$ by $x'$ and enumerate the members of $x'$ in increasing order by $\langle x_i:i<\beta\rangle$. For every $j\in\omega$ let $a_j\in[\beta]^\omega$ be such that $y_j=\{x_i:i\in a_j\}$. Finally, choose $w\in [A]^{\aleph_0\text{-bd}}$ such that ${\rm sup}(x')<{\rm min}(w), {\rm otp}(w)=\omega$ and $r_{{\rm sup}(w)}(w)=\langle a_j:j\in\omega\rangle$.

We let $x=x'\cup w$. Notice that $x\in [A]^{\aleph_0\text{-bd}}$ and the order type of $x$ is $\beta+\omega$ for some limit ordinal $\beta$. By the very definition of $h$ we have $h(x)=\eta$, so $h``[A]^{\aleph_0\text{-bd}}=\alpha$ as $\eta$ was arbitrary. However, this contradicts the choice of $A$, so we are done.

\hfill \qedref{solreferee}

\newpage

\providecommand{\bysame}{\leavevmode\hbox to3em{\hrulefill}\thinspace}
\providecommand{\MR}{\relax\ifhmode\unskip\space\fi MR }
\providecommand{\MRhref}[2]{%
  \href{http://www.ams.org/mathscinet-getitem?mr=#1}{#2}
}
\providecommand{\href}[2]{#2}

\end{document}